# BARYCENTERS OF MEASURES TRANSPORTED BY STOCHASTIC FLOWS

By Marc Arnaudon and Xue-Mei Li[1]

*Université de Poitiers and Nottingham Trent University*

We investigate the evolution of barycenters of masses transported by stochastic flows. The state spaces under consideration are smooth affine manifolds with certain convexity structure. Under suitable conditions on the flow and on the initial measure, the barycenter $\{Z_t\}$ is shown to be a semimartingale and is described by a stochastic differential equation. For the hyperbolic space the barycenter of two independent Brownian particles is a martingale and its conditional law converges to that of a Brownian motion on the limiting geodesic. On the other hand for a large family of discrete measures on suitable Cartan–Hadamard manifolds, the barycenter of the measure carried by an unstable Brownian flow converges to the Busemann barycenter of the limiting measure.

**1. Introduction.** We consider the motion of a mass moving according to the law of a random flow. This can be used to model the motion of passive tracers in a fluid, for example, the spread of oil spilled in an ocean. Such motion can be assumed to obey a stochastic flow where particles at nearby points are correlated, for example, isotropic flows as investigated in [22, 27] or a general semimartingale flow as considered in [5]. The evolution of pollution clouds in the atmosphere or a gas of independent particles, on the other hand, can be described as blocks of masses moving according to the laws of independent stochastic flows. Here we propose to study the dynamics of masses transported by stochastic flows by investigating the motion of their centers of mass. As the medium in which the liquid travels is not necessarily homogeneous or flat, it makes sense to work on a nonlinear space, for example, on a manifold diffeomorphic to $\mathbb{R}^n$ but with different

Received April 2003; revised March 2004.
[1]Supported by NSF Grant DMS-00-72387.

geometric structure. All measures considered here are normalized to have mass 1.

The center of mass for a measure on $\mathbb{R}^n$ is the minimizer of the square distance function averaged with respect to the measure. This definition is used to define martingales on $\mathbb{R}^n$. The same minimizing problem can be considered on Riemannian manifolds. Although the distance in $L^2$ is the traditional object to minimize, we can also make sense of finding the minimizer of the distance function in $L^1$. Busemann barycenters, for example, of measures on the boundaries of hyperbolic spaces are minimizers of Busemann distance functions in $L^1$. The Busemann barycenter was investigated for continuous measures for its application in analysis (see, e.g., [3]), while we are interested in both discrete and continuous measures. We shall relate the center of a mass pushed by a stochastic flow, in the limit, to the Busemann barycenter of its limiting measure.

From the point of view of minimizing assumptions on the state space, we note below that a center of mass can be defined using the concept of geodesics. For this we only need a linear connection $\nabla$ on a smooth manifold, and the pair $(M, \nabla)$ shall be called an *affine manifold*. Denote by $\gamma_u$ the geodesic with initial velocity $u \in T_xM$. Denote by $\exp: TM \to M$ the exponential map $\exp_x u = \gamma_u(1)$, if it exists.

A point $x \in M$ is an *exponential barycenter* of $\mu$ if it is a solution to

$$(1.1) \qquad \int_M \exp_x^{-1} y \, \mu(dy) = 0.$$

In the sequel, by barycenter we always mean exponential barycenter. Note that (1.1) always makes sense locally and therefore for measures with sufficiently small support. However, to consider more general measures we are obliged to work on convex manifolds. An affine manifold $(M, \nabla)$ is said to be *convex* if for every pair of points $x, y \in M$ there exists a unique geodesic, defined using $\nabla$, joining $x$ and $y$ and this geodesic depends smoothly on $x$ and $y$.

*Main results.* In Section 2 the existence and uniqueness of exponential barycenters of a probability measure with compact support are established for CSLCG manifolds. A CSLCG manifold is an affine manifold satisfying certain convexity conditions. Furthermore, a stochastic differential equation governing the motion of the barycenter is given.

To study the large time behavior of the barycenter, we first investigate the toy model of the empirical measure $\mu_t = \frac{1}{n} \sum_{i=1}^n \delta_{X_t^i}$ of $n$ independent Brownian particles $(X_t^i)$ on an hyperbolic space $\mathbb{H}$. For $n = 2$ this is explicitly studied in Section 3 and for $n \geq 3$ this can be considered as a special case of the discussion in Section 4 and can be proved analogously. For two particles, when $t$ gets large the barycenter gets close to a Brownian motion of variance $1/2$ on the geodesic with asymptotic directions $X_\infty^1$ and $X_\infty^2$, where $X_\infty^i$ is the limit of $X_t^i$. For more than two particles, the barycenter converges to a unique point in the manifold. This result seems to be surprising. However, it has a link to the elementary Steiner problem of finding a point which minimizes the sum of its distances to $n$ given vertices: there is a unique minimizer if and only if the vertices do not all lie on the same line. If we replace the usual distance in this problem by the Busemann distance function and look for the point which minimizes the sum of the distances to the $n$ given points on the boundary, we shall obtain the Busemann barycenter of the empirical measure of the limiting points. The unique point the barycenter of $\mu_t$ converges to is this Busemann barycenter. For $n = 2$, the two limiting points can be joined by a geodesic and every point on the geodesic is a minimizer for the Busemann distance function.

In Section 4 the state space under consideration is a Cartan–Hadamard manifold with pinched negative curvature on which the Brownian motions satisfy a law of large numbers. In this case a Brownian flow $F_t(y)$ converges almost surely for every starting point $y$ to a process $F_\infty(y)$ on the visibility compactification of $M$. Furthermore, under suitable conditions, the limiting process separates the initial points and the uniform $o(t)$ fluctuations for the Brownian motions vanish in the limit. Thus the $L^2$ minimizer of the distance function is the same as the $L^1$ minimizer, the Busemann barycenter, of the limiting measure on the boundary. More precisely, suppose $\mu_0$ has finite support and all the weights are smaller than $1/2$. Then the barycenter of $\mu_t$ converges almost surely to a random variable $Z_\infty \in M$, the Busemann barycenter of $F_\infty(\mu_0)$, characterized by

$$\lim_{t\to\infty} \int_M \dot\varphi(Z_\infty, F_t(y))(0) \, d\mu(y) = 0,$$

where $(\varphi(z,y)(s), s \geq 0)$ is the unit speed geodesic connecting $z$ and $y$ starting at $z$ and $\dot\varphi(z,y)(0) = \frac{d}{ds}|_{s=0}\varphi(z,y)(s)$.

*Notation.* Throughout the paper, a standard filtered probability space $(\Omega, \mathcal{F}, (\mathcal{F}_t)_{t\geq 0}, \mathbb{P})$ is fixed. If $(X_t)$ is an $M$-valued continuous semimartingale and $\alpha$ is a $C^1$ section of $T^*M$, we denote by $\int_0^t \langle \alpha(X_r), \delta X_r \rangle$ and $\int_0^t \langle \alpha(X_r), d^\nabla X_r \rangle$, respectively, the Stratonovich and Itô integrals of $\alpha$ along $X$. The second integral requires a connection $\nabla$ on $M$. If $X_t(\omega)$ takes its values in the domain of a local chart for any $t \in [S(\omega), T(\omega)[$ where $S$ and $T$ are stopping times, then

$$\int_S^T \langle \alpha(X_t), \delta X_t \rangle = \int_S^T \alpha_i(X_t) \, dX_t^i + \frac{1}{2} \int_S^T \frac{\partial \alpha_i(X_t)}{\partial x^j} \, d\langle X^i, X^j \rangle_t$$

and

$$\int_S^T \langle \alpha(X_t), d^\nabla X_t \rangle = \int_S^T \alpha_i(X_t)(dX_t^i + \tfrac{1}{2}\Gamma_{jk}^i(X_t) \, d\langle X^j, X^k \rangle_t),$$

where $\Gamma^i_{jk}$ are the Christoffel symbols of $\nabla$. The Itô integral can be defined via the Stratonovich integral as follows: let $//_{0,t}: T_{X_0}M \to T_{X_t}M$ be the stochastic parallel transport along $(X_t)$ and let $Z_t = \int_0^t //_{0,r}^{-1} \delta X_r$ be the anti-development of $(X_t)$, a $T_{X_0}M$-valued semimartingale. Then

$$\int_0^t \langle \alpha(X_s), d^\nabla X_s \rangle = \int_0^t \langle \alpha(X_s) \circ //_{0,s}, dZ_s \rangle.$$

A semimartingale $(X_t)$ is called a $\nabla$-*martingale* if and only if, for every smooth section $\alpha$ of $T^*M$, $\int_0^t \langle \alpha(X_s), d^\nabla X_s \rangle$ is a local martingale. Equivalently the stochastic anti-development of $(X_t)$ by $\nabla$ is a local martingale. The reader may consult [7, 8, 9] and [12] for general reference.

## 2. Exponential barycenters of measures transported by flows.

2.1. *Preliminaries.* If the geodesics determined by two linear connections are the same, then the two connections differ by a tensor of type $(1,2)$. In particular, given any connection, we can always subtract from it half of its torsion $T$ to obtain a torsion-free connection with the same geodesics. Since convexity and barycenters are determined by geodesics, we assume that $\nabla$ is a torsion-free connection. We do not impose the condition that our connection be metric with respect to any Riemannian metric on $M$. However, for calculations we shall fix an auxiliary Riemannian metric on $M$ with the corresponding Riemannian distance function denoted by $\rho$. In the sequel when the term "*Riemannian manifold*" is used it is considered as an affine manifold with respect to the Levi–Civita connection. A smooth local diffeomorphism between two manifolds $M$ and $(\tilde M, \tilde\nabla)$ induces a connection $\nabla$ on $M$. An affine map $\phi$ between two affine manifolds $(M, \nabla)$ and $(\tilde M, \tilde\nabla)$ is a smooth map such that for all smooth functions $f: \tilde M \to \mathbb{R}$, $\nabla d(f \circ \phi) = \phi^*(\tilde\nabla df)$. An affine map on $M$ preserves geodesics. We shall frequently use the following:

1. A point $x$ in $M$ is a barycenter of $\mu$ if and only if

(2.1) $$H(x) \equiv H_\mu(x) = \int_M \dot\gamma(x,y)(0)\mu(dy) = 0$$

   for geodesics $(\gamma(x,y)(s): 0 \le s \le 1)$ connecting $x$ and $y$.
2. If $z$ is the barycenter of $\mu$ and $G: M \to N$ an affine map between two convex affine manifolds, then $G(z)$ is the barycenter of the push forward measure $G(\mu)$ by $G$.
3. For a smooth Riemannian manifold let $\mu$ be a probability measure on $M$ such that $f(x) = \int_M \text{dist}^2(x,y)\,\mu(dy)$ is finite for some (hence for all) $x$ in $M$. Then the exponential barycenters of $\mu$ are the critical points of the function $f$ ([10], Proposition 3). Moreover, if $x$ is an exponential barycenter of $\mu$ and $h$ is a convex function on $M$, then we have Jensen's-type inequality ([10], Proposition 2)

$$h(x) \leq \int_M h(y)\,\mu(dy). \tag{2.2}$$

A convex manifold $M$ is diffeomorphic to an open set of $\mathbb{R}^m$, $m = \dim(M)$, and there exists a neighborhood $U$ of the null section in $TM$ such that

$$(x,u) \in U \mapsto (x, \exp_x u) \in M \times M$$

is a diffeomorphism. Standard convex complete Riemannian manifolds include Cartan–Hadamard manifolds (complete, simply connected manifolds with sectional curvature less than or equal to 0). Examples of incomplete convex manifolds are given by: small geodesic balls in Riemannian manifolds or small balls centered at the origin in an exponential chart in an affine manifold. [Note that every point $x$ in a Riemannian manifold $M$ has a convex geodesic ball $B_r(x)$; in the case of a sphere of radius $r$, a maximal convex set is given by an open half sphere, and its convexity radius is $\pi r/2$.]

Recall that a function $\phi\colon M \to \mathbb{R}$ on an affine manifold $(M, \nabla)$ is *convex* if $\phi \circ \gamma$ is a convex function for all geodesics $\gamma$, that is, for all geodesics $(\gamma(x,y)(t), 0 \leq t \leq 1)$ connecting two different points $x$ and $y$,

$$\phi(\gamma(x,y)(t)) \leq t\phi(x) + (1-t)\phi(y). \tag{2.3}$$

For a $C^2$ function $\phi$, it is convex if and only if $\frac{D}{dt}\frac{d}{dt}\phi(\gamma_t) \geq 0$, where $\frac{D}{dt}$ denotes covariant differentiation with respect to $t$.

For the uniqueness of exponential barycenter we often use a convex function on $M \times M$ which separates points. The definition below is partly borrowed from [9].

DEFINITION 2.1. Consider an affine manifold $(M, \nabla)$ and $M \times M$ with product connection. A convex function $\phi\colon M \times M \to \mathbb{R}_+$ vanishing exactly on the diagonal $\Delta$ of $M \times M$ is called a *separating function* on $M$. Let $p \in 2\mathbb{N}$. A manifold which carries a smooth separating function $\phi$ such that

$$c\rho^p \leq \phi \leq C\rho^p \tag{2.4}$$

for some constants $0 < c < C$ and some Riemannian distance function $\rho$ is called a *manifold with p-convex geometry*.

For Cartan–Hadamard manifolds, the Riemannian distance function is a separating function. Note that the square of the distance function is also a separating function and it is smooth. In more general Riemannian manifolds, sufficiently small geodesic balls have 2-convex geometry. For instance, any geodesic ball strictly smaller than an open hemisphere has $p$-convex geometry for some $p$ depending on the radius, as proved by Kendall in [15], and any point in an affine manifold has a convex neighborhood with 2-convex geometry. However, an affine manifold carrying a separating function is not necessarily convex. For example, $\mathbb{R}^m \setminus \{0\}$ is not convex even though the distance function is a separating function. It is more difficult to find convex manifolds which do not carry separating functions; see [16].

Let $\varphi$ be a $C^2$ convex function with $\{\varphi < 0\}$ a relatively compact set. Then any probability measure $\mu$ on $M$ with support included in $\{\varphi < 0\}$ has an exponential barycenter, see [10, 14]. Also if $M$ is an affine manifold with a bounded separating function, then any probability measure on $M$ has at most one exponential barycenter [10, 14].

As a consequence: let $M$ be a convex manifold with convex geometry, and let $\phi$ be a $C^2$ separating function with $\phi_x^{-1}([0, a[)$ relatively compact in $M$ for some fixed $x \in M$ and some $a > 0$ where $\phi_x = \phi(x, \cdot)$ (which is not the case in general). Then any measure supported on the set $\phi_x^{-1}([0, a[)$ has a unique exponential barycenter.

Observe that an open convex subset of a convex manifold with $p$-convex geometry is a manifold with $p$-convex geometry.

2.2. *Existence and uniqueness.*

DEFINITION 2.2. A convex affine manifold $(M, \nabla)$ is said to be CSLCG (convex, with semilocal convex geometry) if every compact subset $K$ of $M$ has a relatively compact convex neighborhood $U_K$ which has $p$-convex geometry for some $p \in 2\mathbb{N}$ depending on $K$.

Equivalently, a convex affine manifold $(M, \nabla)$ is CSLCG if there exists an increasing sequence $(U_n)_{n \geq 1}$ of relatively compact open convex subsets of $M$ such that $M = \bigcup_{n \geq 1} U_n$, and for every $n \geq 1$, $U_n$ has $p$-convex geometry for some $p \in 2\mathbb{N}$ depending on $n$.

This definition is motivated by Propositions 2.4 and 2.7 which state that any probability measure in a CSLCG manifold, of compact support, has a unique barycenter and that in a certain sense, the exponential barycenter of $\mu$ is differentiable as a function of $\mu$.

Note $p$ increases with $K$. Examples of CSLCG manifolds are open hemispheres endowed with the Levi–Civita connection (which do not have $p$-convex geometry for any $p \in 2\mathbb{N}$; see [15]). Examples of geodesically complete CSLCG Riemannian manifolds should be given by manifolds with a pole under curvature conditions to be determined. We conjecture that if $(M, \nabla)$ is a CSLCG manifold, then $(TM, \nabla^c)$ is a CSLCG manifold, where $\nabla^c$ is the complete lift of $\nabla$. Note that on every relatively compact convex subset of $TM$ one can construct a continuous separating function: by [17], uniqueness of martingales with prescribed terminal values implies the existence of such a separating function, and by [1], uniqueness holds. So the conjecture concerns the smoothness of the separating function, and the fact that it satisfies (2.4).

One can find in [16] an example of convex manifold which has not semilocal convex geometry. In the manifold constructed there, there exists a probability measure carried by three points, which possesses four exponential barycenters. Clearly there is no convex neighborhood of these seven points in which we can define a separating function.

LEMMA 2.3. *Let $(M, \nabla)$ be a CSLCG manifold. Every compact convex subset $K$ of $M$ has a convex neighborhood $U$ with a nonnegative $C^1$ convex function $\phi_K$ such that $\phi_K^{-1}(\{0\}) = K$.*

PROOF. Let $U'$ be an open convex relatively compact neighborhood of $K$ with $p$-convex geometry for some $p \in 2\mathbb{N}$. Let $\phi$ be a smooth separating function on $U'$ satisfying $c\rho^p \leq \phi \leq C\rho^p$ where $0 < c < C$, and $\rho$ is the Euclidean distance induced by a global chart, for example, an exponential chart. Define for $x \in U'$

$$\phi_K(x) = \inf\{\phi(y, x), y \in K\}. \tag{2.5}$$

Clearly $\phi_K^{-1}(\{0\}) = K$, since $K$ is compact.

1. First we show that $\phi_K$ is convex on $U'$. Take $(x_1, x_2) \in U' \times U'$. Let $p(x_i)$ be points in the compact set $K$ achieving the minimum: $\phi_K(x_i) = \phi(p(x_i), x_i)$, $i = 1, 2$. Since $\phi$ is convex, for every $t \in [0, 1]$,

$$\phi(\gamma(p(x_1), p(x_2))(t), \gamma(x_1, x_2)(t)) \leq (1-t)\phi(p(x_1), x_1) + t\phi(p(x_2), x_2).$$

On the other hand, by the definition of $\phi_K$, we have since $\gamma(p(x_1), p(x_2))(t) \in K$

$$\phi_K(\gamma(x_1, x_2)(t)) \leq \phi(\gamma(p(x_1), p(x_2))(t), \gamma(x_1, x_2)(t)).$$

Putting these equations together we obtain

$$\phi_K(\gamma(x_1, x_2)(t)) \leq (1-t)\phi_K(x_1) + t\phi_K(x_2),$$

which proves the convexity of $\phi_K$.

2. Next we show that $\phi_K$ is $C^1$ on some convex neighborhood $U$ of $K$ included in $U'$. We first prove that there exists a neighborhood $U \subset U'$ of $K$ on which there is a unique point $p(x)$ such that $\phi_K(x) = \phi(p(x), x)$.

Suppose for $x \in U'$ and $y_0, y_1 \in K$, $\phi_K(x) = \phi(y_0, x) = \phi(y_1, x)$. Set $y(t) = \gamma(y_0, y_1)(t)$. Then necessarily, for every $t \in [0, 1]$, $\phi(y(t), x) = \phi_K(x)$.

By the Hadamard lemma (see, e.g., [4], Corollaire 3.1.9 and Exercice 3.1.10) we can write in the global chart

$$\phi(y, x) = \sum_{i_1, \ldots, i_p = 1}^{m} a_{i_1 \ldots i_p}(x) \prod_{j=1}^{p} (y^{i_j} - x^{i_j})$$

$$+ \sum_{i_1,\ldots,i_{p+1}=1}^{m} b_{i_1\ldots i_{p+1}}(y,x) \prod_{j=1}^{p+1}(y^{i_j}-x^{i_j}),$$

where $a_{i_1\ldots i_p}$ and $b_{i_1\ldots i_{p+1}}$ are smooth functions.

Write $f(t)=\phi(y(t),x)$. It is a constant function. On the other hand, we can differentiate $p$ times the function $f$ with the expression of $\phi$ in the chart. Using

$$\ddot{y}^i(t)=-\sum_{j,k=1}^{m}\Gamma^i_{jk}(y(t))\dot{y}^j(t)\dot{y}^k(t),$$

where $\Gamma^i_{jk}$ are the Christoffel symbols of $\nabla$ in the chart, we see that $f^{(p)}$ is of the following format:

$$f^{(p)}(t)=\sum_{i_1,\ldots,i_p=1}^{m} a_{i_1\ldots i_p}(x)\prod_{j=1}^{p}\dot{y}^{i_j}(t)+g(y(t),\dot{y}(t),x),$$

where $g$ is a smooth function. Since $U'$ is a relatively compact convex open neighborhood of $K$, we see by the Hadamard lemma and explicit calculation that

$$|g(y,z,x)|\leq C'\|y-x\|\|z\|^p$$

for some positive constant $C'$. Now since $c\rho^p\leq\phi$, we have for every $z\in\mathbb{R}^m$

$$\sum_{i_1,\ldots,i_p=1}^{m} a_{i_1\ldots i_p}(x)\prod_{j=1}^{p}z^{i_j}\geq c\|z\|^p.$$

We choose $U$ of the form $U=\{x'\in U',\phi_K(x')<\varepsilon\}$, with $\varepsilon\in\,]0,c^{p+1}/(C')^p[$, such that $U$ is convex and relatively compact in $U'$. Consequently for $x\in U$ and $y\in K$, $\|y-x\|^p\leq\phi(x,y)/c\leq\varepsilon/c$ and

$$f^{(p)}(t)\geq c\|\dot{y}(t)\|^p-C'\|y(t)-x\|\cdot\|\dot{y}(t)\|^p$$
$$\geq(c-(\varepsilon/c)^{1/p}C')\|\dot{y}(t)\|^p\qquad\forall\,t\in[0,1],$$

where $c-(\varepsilon/c)^{1/p}C'>0$. But $f$ is constant, so we have $\dot{y}(t)\equiv 0$ and consequently $y_0=y_1$.

3. For $x\in U$, we let $p(x)$ be the point in $K$ such that $\phi_K(x)=\phi(p(x),x)$. We prove that $p$ is continuous on $U$. If not, let $(x_n)_{n\in\mathbb{N}}$ be a convergent sequence in $U$ with limit $x$ such that $p(x_n)$ does not converge to $p(x)$. Since $K$ is compact, by choosing a subsequence if necessary, we can assume $p(x_n)\to y\in K\setminus\{p(x)\}$. By the continuity of $\phi$, $\phi_K(x_n)\equiv\phi(p(x_n),x_n)\to\phi(y,x)$. On the other hand, since $\phi_K$ is convex on the open set $U$, it is continuous (see, e.g., [11], Proposition 1) and so $\lim\phi_K(x_n)=\phi_K(x)$. Consequently, $\phi(y,x)=\phi_K(x)$ and by the uniqueness given in (2) we obtain $y=p(x)$, a contradiction. So the map $p$ is continuous.

4. We are left to prove that $\phi_K$ is $C^1$. Denote by $d_g\phi_K$ the Gâteaux-differential of $\phi$. Recall that for every $x \in U$ and $v \in T_xM$,

$$d_g\phi_K(x)(v) = \lim_{t\downarrow 0} \frac{\phi_K(\exp tv) - \phi_K(x)}{t}$$

and that $d_g\phi_K(x)$ is convex on $T_xM$ (again by [11], Proposition 1). Let $x \in U$ and $v \in T_xM$. If $t > 0$ is sufficiently small, we have

$$\phi_K(\exp tv) - \phi_K(x) \leq \phi(p(x), \exp tv) - \phi(p(x), x)$$

which yields

$$d_g\phi_K(x)(v) \leq d\phi_{p(x)}(x)(v),$$

$\phi_y$ denoting the map $\phi(y, \cdot)$. Since $d\phi_{p(x)}(x)$ is linear, this inequality implies $d_g\phi_K(x) = d\phi_{p(x)}(x)$ on $T_xM$. Otherwise we would have $d_g\phi_K(x)(v) < d\phi_{p(x)}(x)(v)$ for some $v \in T_xM$, which would give by convexity of $d_g\phi_K(x)$

$$d_g\phi_K(x)(-v) \geq -d_g\phi_K(x)(v) > -d\phi_{p(x)}(x)(v) = d\phi_{p(x)}(x)(-v),$$

a contradiction. So $d_g\phi_K(x) = d\phi_{p(x)}(x)$ on $T_xM$. The differentiability of $\phi_K$ then easily comes from the inequalities, for all $v \in T_xM$ sufficiently close to $0$,

$$0 \leq \phi_K(\exp v) - \phi_K(x) - d\phi_{p(x)}(x)(v)$$
$$\leq \phi_{p(x)}(\exp v) - \phi_{p(x)}(x) - d\phi_{p(x)}(x)(v).$$

That $\phi_K$ is $C^1$ comes from the continuity of $p$. □

PROPOSITION 2.4. *Let $M$ be a CSLCG manifold. Every probability measure $\mu$ on $M$ with compact support has a unique exponential barycenter.*

PROOF. The proof is a slight modification of that of Proposition 5 in [10]. For uniqueness let $x$ and $x'$ be two exponential barycenters of $\mu$. Set $K = \text{supp}(\mu) \cup \{x, x'\}$. Let $\phi$ be a separating function defined on a convex relatively compact neighborhood of $K$. Writing $\nu$ the measure $\mu$ pushed forward by $y \mapsto (y, y)$, $(x, x')$ is an exponential barycenter of $\nu$, so

$$\rho^p(x, x') \leq \frac{1}{c}\phi(x, x') \leq \frac{1}{c}\int \phi(y, y')\, d\nu(y, y') = \frac{1}{c}\int \phi(y, y)\, d\mu(y) = 0$$

giving $x = x'$.

For the existence let $K$ be a convex compact subset of $M$ containing the support of $\mu$ (we know that this exists in a CSLCG manifold). By Lemma 2.3 there exists a $C^1$ convex nonnegative function $\phi_K$ defined on a relatively compact open neighborhood $U$ of $K$, such that $\phi_K^{-1}(\{0\}) = K$. Let $\varepsilon > 0$ satisfy $\phi_K^{-1}([0, \varepsilon[) \subset U$. We apply [10], Proposition 5, to the function $\phi_K - \varepsilon/2$ to see that $\mu$ has an exponential barycenter in $\phi_K^{-1}([0, \varepsilon/2[)$ (note [10], Proposition 5 is still valid with a $C^1$ convex function instead of a $C^2$ convex function). □

2.3. *Mass moved by a smooth vector field.* Let $(M, \nabla)$ be a smooth convex affine manifold with $\nabla$ torsion free. As before, for $x$, $y$ in $M$ let $(\gamma(x,y)(s), 0 \leq s \leq 1)$ be the geodesic, with respect to $\nabla$, with end points $x$ and $y$. The geodesic $\gamma(x,y)$ shall be abbreviated as $\gamma$ where there is no risk of confusion. Set

$$J(u,v)(s) = T(\gamma(\cdot,\cdot)(s))(u,v), \qquad u \in T_x M,\ v \in T_y M,\ s \in \mathbb{R},$$

where $T\gamma(s)(u,v)$ denotes the derivative of $\gamma(s)$ in the direction of $(u,v)$. The map $(J(u,v)(s), 0 \leq s \leq 1)$ is the Jacobi field satisfying boundary condition $J(u,v)(0) = u$ and $J(u,v)(1) = v$. Write $\dot{J}(u,v)(s) = (\nabla_{d/ds} J(u,v))(s)$. For every $x, y \in M$, we define the linear map

$$(2.6) \qquad \psi_{(x,y)} : T_x M \to T_x M, \qquad u \mapsto \dot{J}(u, 0_y)(0),$$

where $0_y$ is the zero tangent vector in $T_y M$.

On $TM$ one can define a connection $\nabla^c$, called the complete lift of $\nabla$ (see, e.g., [26]). It is a torsion-free connection since $\nabla$ is, so it is characterized by its geodesics, which are the Jacobi fields of the connection $\nabla$. The canonical projection $\pi : (TM, \nabla^c) \to (M, \nabla)$ is affine. We easily see that $(TM, \nabla^c)$ is convex if $(M, \nabla)$ is and the only geodesic joining $u$ and $v$ in $TM$ is $(J(u,v)(s), 0 \leq s \leq 1)$. Given a $C^1$ vector field $A$, consider the map $A : M \to TM$. It induces a measure $A(\mu)$ on $TM$. A point $v \in TM$ is a barycenter of a measure $A(\mu)$ if and only if $\pi(v)$ is a barycenter of $\mu$ and

$$(2.7) \qquad \int_M \dot{J}(v, A(x))(0)\, d\mu(x) = 0.$$

A measure $\mu$ is differentiable along a vector field $A$ if its pushed forward measure $S_t(\mu)$ is weakly differentiable in $t$, where for $x \in M$, $S_t(x)$ is the integral curve of $A$ starting from $x$. The following lemma asserts that the map $\mu \to b(\mu)$, $b(\mu)$ the barycenter of $\mu$, is a differentiable map in the above sense.

PROPOSITION 2.5. *Suppose $(M, \nabla)$ is a CSLCG manifold. For $i = 1, \ldots, n$, let $\mu^i$ be a probability measure on $M$ with compact support $K$, let $A^i$ be a smooth $C^1$ vector field and let $(S^i_t(x), 0 \leq t < \tau)$ be its integral curve starting from $x$. For $p_i \geq 0$ with $\sum_{i=1}^n p_i = 1$, set $\mu_t = \sum_{i=1}^n p_i S^i_t(\mu^i)$. If $z(t)$ is the exponential barycenter of $\mu_t$:*

(a) *then $z(-)$ is differentiable at $t = 0$ and $\dot{z}(0)$ is the exponential barycenter of $\sum_{i=1}^n p_i A^i(\mu^i)$ with respect to the complete lift $\nabla^c$ of $\nabla$;*

(b) *there is $C_K > 0$ depending only on $K$ and the arbitrary metric such that*

$$\|\dot{z}(0)\| \leq C_K \sup_{\substack{x \in K \\ 1 \leq i \leq n}} \|A^i(x)\|.$$

PROOF. Let $M'$ be a convex relatively compact open neighborhood of $K$ with $p$-convex geometry ($p \in 2\mathbb{N}$). Let $\phi$ be a smooth separating function on $M'$ satisfying $c\rho^p \leq \phi \leq C\rho^p$ where $\rho$ is a Riemannian distance function, $0 < c < C$. For $t < \tau$ the barycenter $z(t)$ exists and by the definition of barycenters, $(z(0), z(t))$ is the exponential barycenter of $\sum_{i=1}^n p_i(S_0^i, S_t^i)(\mu^i)$ with respect to the product connection.

Apply (2.2) to see $\phi(z(0), z(t)) \leq \sum_{i=1}^n p_i \int_K \phi(S_0^i(x), S_t^i(x)) \, d\mu^i(x)$, which implies

$$
\begin{aligned}
\rho(z(0), z(t)) &\leq \left[ \frac{1}{c} \phi(z(0), z(t)) \right]^{1/p} \\
&\leq \left( \frac{C}{c} \right)^{1/p} \left( \sum_{i=1}^n p_i \int_K \rho^p(S_0^i(x), S_t^i(x)) \, d\mu^i(x) \right)^{1/p} \\
&\leq \left( \frac{C}{c} \right)^{1/p} t \sup_{i \in \{1,\ldots,n\}} \sup_{\substack{x \in K \\ 0 \leq r \leq t}} \|A^i(S_r^i(x))\|,
\end{aligned}
\tag{2.8}
$$

where $\|\cdot\|$ is the Riemannian metric corresponding to $\rho$. Since all Riemannian metrics on $M'$ are comparable,

$$
\begin{aligned}
\left\| \frac{1}{t} \exp_{z(0)}^{-1} z(t) \right\| &\leq \frac{c_0}{t} \rho(z(0), z(t)) \\
&\leq c_0 \left( \frac{C}{c} \right)^{1/p} \sup_{i \in \{1,\ldots,n\}} \sup_{\substack{x \in K \\ 0 \leq r \leq t}} \|A^i(S_r^i(x))\|,
\end{aligned}
\tag{2.9}
$$

some constant $c_0$ and so the family $\{\frac{1}{t} \exp_{z(0)}^{-1} z(t), t < \tau\}$ is bounded and thus has a limit point which shall be denoted by $u$. Let $t_k$ be a sequence going to 0 such that $\lim_{k \to \infty} \frac{1}{t_k} \exp_{z(0)}^{-1} z(t_k) = u$. Each $z_{t_k}$ satisfies

$$
0_{z(t_k)} = \sum_{i=1}^n p_i \int_M \dot{\gamma}(z(t_k), S_{t_k}^i(x))(0) \, d\mu^i(x)
$$

and the covariant derivative of $\dot{\gamma}(z(t), S_t^i(x))(0)$ along the subsequence $t_k$ satisfies

$$
\left. \frac{D}{dt_k} \right|_{t_k \to 0} \dot{\gamma}(z(t_k), S_{t_k}^i(x))(0) = \dot{J}(u, A^i(x))(0).
$$

Integrate the identity over $M$ and use the dominated convergence theorem to see

$$
0_{z(0)} = \sum_{i=1}^n p_i \int_M \dot{J}(u, A^i(x))(0) \, d\mu^i(x)
$$

which together with

$$0_{z(0)} = \sum_{i=1}^n p_i \int_M \dot\gamma(z(0),x)(0)\,d\mu^i(x)$$

shows that $u$ is an exponential barycenter of $\sum_{i=1}^n p_i A^i(\mu^i)$ with respect to $\nabla^c$.

Let $\pi: TM \to M$ be the canonical projection. To prove the uniqueness consider the map

$$T^{\otimes p}\phi: TM' \to \mathbb{R}_+, \qquad v \mapsto \lim_{t\to 0^+} \frac{1}{t^p}\phi(\pi(v), \exp tv),$$

which is $\nabla^c$-convex and $T^{\otimes p}\phi(v) \leq C\|v\|^p$. Furthermore the map

$$(v,v') \mapsto T^{\otimes p}\phi(v'-v)$$

is convex on $E := \{(v,v') \in TM' \times TM',\ \pi(v) = \pi(v')\}$ for the connection $\nabla^c \otimes \nabla^c$ (note $E$ is a totally geodesic submanifold), and vanishes exactly on the diagonal of $TM' \times TM'$. See [1] for the proof and details.

Let $u'$ be another exponential barycenter of $A(\mu)$ with respect to $\nabla^c$. Since $\mu$ has a unique barycenter and $\pi$ is affine, $\pi(u') = \pi(u)$. So it makes sense to estimate $T^{\otimes p}\phi(u'-u)$. Set $\nu = \sum_{i=1}^n p_i(A^i, A^i)(\mu^i)$. Then $(u,u')$ is an exponential barycenter for $\nu$, so

$$\begin{aligned}
T^{\otimes p}\phi(u-u') &\leq \int_E T^{\otimes p}\phi(v-v')\,d\nu(v,v') \\
&= \sum_{i=1}^n p_i \int_M T^{\otimes p}\phi(A^i(x) - A^i(x))\,d\mu^i(x) = 0.
\end{aligned}$$

This implies that $u = u'$, which in turn implies that the family $\frac{1}{t}\exp_{z(0)}^{-1} z(t)$ has only one limit point as $t$ goes to $0$ and so it converges to $u$. Thus $t \mapsto z(t)$ is differentiable at $t=0$ with derivative the exponential barycenter of $\sum_{i=1}^n p_i A^i(\mu^i)$ with respect to $\nabla^c$. Finally, using (2.8) and (2.9), letting $t$ go to 0, invoking the compactness of $K$ and the continuity of $S^i$, we get

$$\|\dot z(0)\| \leq C_K \sup_{\substack{i\in\{1,\ldots,n\}\\ x\in K}} \|A^i(x)\| \qquad \text{where } C_K = c_0\left(\frac{C}{c}\right)^{1/p}. \qquad \square$$

The map $\psi_{x,y}$ defined in (2.6) is not invertible in general, as can be seen with two points $x$ and $y$ at a distance $\pi/2$ in a hemisphere. However, we have the following:

LEMMA 2.6. *Let $M$ be a CSLCG manifold, let $K \subset M$ be a compact set and let $\rho$ be a Riemannian metric on $M$. Then there exists a constant* $C'_K > 0$ such that for any probability measure $\mu$ in $M$ with support included in $K$, letting $z$ be the exponential barycenter of $\mu$, the linear map on $T_zM$:

$$u \mapsto \int_M \psi_{(z,x)}(u) \, d\mu(x)$$

is invertible with inverse bounded by $C'_K$ (with respect to the metric $\rho$).

PROOF. Since $M$ is convex, the map $v \mapsto \dot{J}(0_z, v)(0)$ from $T_xM$ to $T_zM$ is injective and so is surjective. Its inverse $-A(x, \cdot)$ is well defined. For each $u \in T_zM$ fixed, $A(\cdot, u)$ is the vector field in $M$ characterized by $\dot{J}(0_z, A(x, u))(0) = -u$. Let $b(u) \in T_zM$ be the exponential barycenter of $A(u)(\mu)$ (by Proposition 2.5, it is well defined). By (2.7)

$$0 = \int_M \dot{J}(b(u), A(x,u))(0)\mu(dx)$$
$$= \int_M \dot{J}(b(u), 0)(0)\mu(dx) + \int_M \dot{J}(0, A(x,u))(0)\mu(dx).$$

By the definition $\psi_{(z,x)}(v) = \dot{J}(v, 0_x)(0)$ for $v \in T_zM$,

$$\left(\int_M \psi_{(z,x)}(\cdot) \, d\mu(x)\right)(b(u)) = -\int_M \dot{J}(0_z, A(x,u))(0) \, d\mu(x) = u.$$

We have $\sup_{x \in K} \|A(x,u)\| \leq C\|u\|$ for some constant $C$ since $K$ is compact. Apply Proposition 2.5 to see

$$(2.10) \qquad \|b(u)\| \leq C_K \sup_{x \in K} \|A(x,u)\| \leq C_K C \|u\|.$$

Thus $\int_M \psi_{(z,x)}(\cdot) \, d\mu(x)$ is invertible with inverse satisfying

$$\left\|\left(\int_M \psi_{(z,x)}(\cdot) \, d\mu(x)\right)^{-1}(u)\right\| \leq \|b(u)\| \leq C'_K\|u\| \qquad \text{for } C'_k = CC_K. \quad \square$$

As a direct consequence of Lemma 2.6 and its proof, we make more precise the result of Proposition 2.5 and give an expression for $\dot{z}(t)$.

PROPOSITION 2.7. *Assume $(M, \nabla)$ is a CSLCG manifold with the connection $\nabla$ torsion free. For $n \geq 1$ and $i = 1, \ldots, n$, let $p_i$ be positive numbers satisfying $\sum_{i=1}^n p_i = 1$ and let $\mu^i$ be probability measures on $M$ with compact support. Let*

$$S^i_\cdot(\cdot) : [0, \infty[ \times \mathrm{supp}(\mu^i) \to M$$

*be a map of class $C^{k,0}$ for some $k \geq 1$ [i.e., $C^k$ in $t$ with derivatives up to order $k$ with respect to $t$ continuous in $(t,x)$]. Then the exponential barycenter $z(t)$ of $\mu_t := \sum_{i=1}^n p_i S^i_t(\mu^i)$ exists and is unique. Furthermore:*

(i) $t \mapsto z(t)$ is $C^k$;
(ii) *the barycenter $z(t)$ is characterized by the identity*

$$\sum_{i=1}^n p_i \int_M \dot J(\dot z(t), 0_{S_t^i(x)})(0)\, d\mu^i(x) = -\sum_{i=1}^n p_i \int_M \dot J(0_{z(t)}, \dot S_t^i(x))(0)\, d\mu^i(x)$$

*where $0_z$ is the zero vector in $T_zM$;*
(iii) *$z(t)$ solves the ordinary differential equation*

(2.11)
$$\dot z(t) = -\sum_{i=1}^n p_i \int_M \left(\sum_{j=1}^n p_j \int_M \psi_{(z(t)S_t^j(x))}\, d\mu^j(x)\right)^{-1} \\ \times \dot J(0_{z(t)}, \dot S_t^i(y))(0)\, d\mu^i(y),$$

*where $\psi$ is defined by* (2.6).

Taking $\mu^i$ to be Dirac distributions, we get from Proposition 2.7 the following:

PROPOSITION 2.8. *If $(M, \nabla)$ is a CSLCG manifold with the connection $\nabla$ torsion free let $n \geq 1$ and $p = (p_1, \ldots, p_n)$ be an $n$-tuple of positive numbers satisfying $\sum_{i=1}^n p_i = 1$. Then the map $G_p$ defined by*

$$G_p: M^n \to M, \qquad (x_1, \ldots, x_n) \mapsto b\left(\sum_{i=1}^n p_i \delta_{x_i}\right),$$

*is smooth with derivative,*

$$TG_p(u_1, \ldots, u_n) = -\sum_{i=1}^n p_i \left(\sum_{j=1}^n p_j \psi_{(G(x_1,\ldots,x_n),x_j)}\right)^{-1} \dot J(0_{G(x_1,\ldots,x_n)}, u_i)(0),$$

*where $(u_1, \ldots, u_n) \in T_{x_1}M \times \cdots \times T_{x_n}M$.*

PROOF. Just note if $S_t^i(x_i) = \exp(tu^i)$, then $\dot z(0) = TG_p(u^1, \ldots, u^n)$. □

We want to generalize Proposition 2.7 to the case where for every $x$ the flow is a semimartingale. We begin with a simple case, which is an immediate corollary of Proposition 2.8 since $Z_t$ as a smooth function of semimartingales is a semimartingale and stochastic calculus applies.

PROPOSITION 2.9. *Let $(M, \nabla)$ be a CSLCG manifold. Let $n \geq 1$, $p = (p_1, \ldots, p_n)$ be an $n$-tuple of positive numbers satisfying $\sum_{i=1}^n p_i = 1$, $x_1 \ldots x_n$* be $n$ points of $M$ and $\mu_0 = \sum_{i=1}^{n} p_i \delta_{x_i}$. For every $i$, let $X_t^i$ be a semimartingale started at $x_i$. Then the exponential barycenter $Z_t$ of $\sum_{i=1}^{n} p_i \delta_{X_t^i}$ is a semimartingale and satisfies

$$\delta Z_t = -\sum_{i=1}^{n} p_i \bigg( \sum_{j=1}^{n} p_j \psi_{(Z_t, X_t^j)} \bigg)^{-1} \dot{J}(0_{Z_t}, \delta X_t^j)(0).$$

If $X_t^i = F_t(x^i)$ where $F_t(x)$ is a semimartingale flow, then

$$\delta Z_t = -\int_M \bigg( \int_M \psi_{(Z_t, F_t(x))} \, d\mu_0(x) \bigg)^{-1} \dot{J}(0_{Z_t}, \delta F_t(y))(0) \, d\mu_0(y).$$

When $\mu_0$ is carried by two points we do not need to assume that the manifold has locally convex geometry: convexity is sufficient. If $\mu_0 = \frac{1}{2}(\delta_{x_0} + \delta_{y_0})$, $X_t = F_t(x_0)$, $Y_t = F_t(y_0)$ for a stochastic flow $F_t$, then $Z_t$ is a semimartingale and

(2.12) $\delta Z_t = -(\psi_{(Z_t, X_t)} + \psi_{(Z_t, Y_t)})^{-1} (\dot{J}(0_{Z_t}, \delta X_t)(0) + \dot{J}(0_{Z_t}, \delta Y_t)(0)).$

2.4. *Mass pushed by a random flow.* Next we consider the case where $F_t$ is a local semimartingale flow of homeomorphisms (which we abbreviate as semimartingale flow). Following Kunita ([18], Section 4.7) with small modifications, see also [20, 21], let $a(t, x, y, \omega)$ be a predictable process with values in the tensor product $T_x M \otimes T_y M$, let $b(t, x, \omega)$ be a predictable process with values in $T_x M$ and let $A_t$ be a continuous adapted real-valued increasing process. We say that $F_t(x)$ is a semimartingale flow with characteristic $(a(t, x, y, \omega), b(t, x, \omega), A_t)$ if

$$dF_t(x) \otimes dF_t(y) = a(t, F_t(x), F_t(y)) \, dA_t \qquad \forall x, y \in M,$$

and the drift of $F_t(x)$ is $b(t, F_t(x)) \, dA_t$, that is, for every 1-form $\alpha \in T^*M$,

$$\int \langle \alpha(F_t(x)), d^\nabla F_t(x) - b(t, F_t(x)) \, dA_t \rangle$$

is a local martingale.

PROPOSITION 2.10. *Let $(M, \nabla)$ be a CSLCG manifold with $\nabla$ torsion free and let $\mu_0$ be a probability measure with compact support. Assume $F_t(x)$ is a semimartingale flow with characteristic $(a(t, x, y, \omega), b(t, x, \omega), A_t)$, such that $a(t, -, -, \omega)$ and $b(t, -, \omega)$ are $C^1$ with derivatives a.s. locally uniformly bounded in time and space. Then the exponential barycenter $Z_t$ of $\mu_t$ exists and is unique. Furthermore:*

(i) *the process $Z_t$ is a semimartingale on $[0, \tau[$, where $\tau$ is an almost surely positive random explosion time;*

(ii) *the barycenter $Z_t$ is characterized by the identity*

$$\int_M \dot{J}(\delta Z_t, 0_{F_t(x)})(0)\, d\mu_0(x) = -\int_M \dot{J}(0_{Z_t}, \delta F_t(x))(0)\, d\mu_0(x);$$

(iii) $Z_t$ *solves the stochastic differential equation*

$$(2.13)\quad \delta Z_t = -\int_M \left(\int_M \psi_{(Z_t, F_t(x))}\, d\mu_0(x)\right)^{-1} \dot{J}(0_{Z_t}, \delta F_t(y))(0)\, d\mu_0(y).$$

REMARK 2.11. The assumptions on the semimartingale flow can be weakened. We can assume, instead, that $(a(t, x, y, \omega), b(t, x, \omega), A_t)$ belongs to the class $B^{0,1}$ in the sense of [18], adapted to a manifold.

PROOF OF PROPOSITION 2.10. Fix an arbitrary Riemannian distance $\rho$ on $M$. For simplicity we choose $a$ and $b$ bounded (by [18], remark on page 85, this amounts to changing $A_t$). Using a change of time, we may assume that $A_t = t$. With the assumptions on the local characteristics of the flow $F_t(x)$, we can choose a version of $F_t(x)$ jointly continuous in $(t, x)$. Let $(U_n)_{n \geq 1}$ be an increasing sequence of relatively compact convex open subsets of $M$ containing $\mathrm{supp}(\mu_0)$ and such that $\bigcup_n U_n = M$. Define

$$\tau_n = \inf\{t > 0, F_t(x) \notin U_n \text{ for some } x \in \mathrm{supp}(\mu_0)\}.$$

Then since the map $(t, x) \mapsto F_t(x)$ is almost surely continuous and $\mu_0$ has compact support, $\tau_n$ is an increasing sequence of stopping times converging to the explosion time $\tau$, which is almost surely positive.

The map $(t, y) \mapsto F_t(y)(\omega)$ is almost surely uniformly continuous on $[0, R] \times \mathrm{supp}(\mu_0)$ for every $R > 0$ smaller than $\tau$. As a consequence, the barycenter $Z_t$ of $\mu_t$ is defined on $[0, \tau[$, and is continuous and adapted.

We shall show that $Z_t$ is a semimartingale and is given precisely by the equations in parts (ii) and (iii) up to each time $\tau_n$. Letting $n$ tend to infinity will prove that the conclusion of the proposition holds.

For the local result we begin with the following lemma which can be considered as an extension of Lemma 2.6.

LEMMA 2.12. *Let $U$ be a relatively compact convex open subset of $M$ containing $\mathrm{supp}(\mu_0)$ and*

$$\tau_U = \inf\{t > 0, F_t(x) \notin U \text{ for some } x \in \mathrm{supp}(\mu_0)\}.$$

*There exists $\varepsilon > 0$ such that if $t \leq \tau_U$ and $\rho(z, Z_t(\omega)) \leq \varepsilon$, then*

$$\int_M \psi_{(z, F_t(x))}\, d\mu_0(x)$$

*is invertible with inverse bounded by a constant depending only on $U$.*

PROOF. Consider the set $(\bar{U})'$ of probabilities on $\bar{U}$, endowed with the weak topology. It is a metrizable space. Since $\bar{U}$ is compact, the map

$$(z, \mu) \mapsto \int_M \psi_{(z,x)} \, d\mu(x)$$

is continuous on $\bar{U} \times (\bar{U})'$ for the product topology. Furthermore, the map from $(\bar{U})'$ to $M$ sending $\mu$ to its exponential barycenter $z(\mu)$ is continuous by Skorokhod theorem: if $\mu_k$ converges to $\mu$, there is a sequence of random variables $X_k$ and a random variable $X$ with values in $\bar{U}$ such that the law of $X_k$ is $\mu_k$, the law of $X$ is $\mu$ and $X_k$ converges almost surely to $X$. Since $\bar{U}$ is compact, $\rho(X_k, X)$ converges to 0 in $L^p$ for every $p \in 2\mathbb{N}$. On the other hand, $(z(\mu_k), z(\mu))$ is the exponential barycenter of the law of $(X_k, X)$ and $\bar{U}$ is $p$-convex for $p$ sufficiently large. So

$$\rho(z(\mu_k), z(\mu)) \leq C \|\rho(X_k, X)\|_p$$

which implies that $z(\mu_k)$ converges to $z(\mu)$. This proves that $\mu \mapsto z(\mu)$ is continuous on $(\bar{U})'$.

Next recall that for every $\mu \in (\bar{U})'$, $\int_M \psi_{(z(\mu),x)} \, d\mu(x)$ is invertible. The image of the compact set $(\bar{U})'$ by the continuous map $\mu \mapsto \int_M \psi_{(z(\mu),x)} \, d\mu(x)$ (w.r.t. the weak topology) is compact. This implies the existence of a neighborhood of the graph of $\mu \mapsto z(\mu)$ of the form $\{(\mu, z) \in ((\bar{U})' \times \bar{U}), \rho(z, z(\mu)) < \varepsilon\}$ for some $\varepsilon > 0$, on which $\int_M \psi_{(z,x)} \, d\mu(x)$ is invertible with inverse in a compact set depending only on $U$. This proves the lemma. $\square$

We continue with the proof of the proposition. For relatively compact set $U_n$, let $\varepsilon_n$ be the constant defined by Lemma 2.12 and set

$$T_n^1 := \tau_n \wedge \inf\left\{t > 0, \rho(Z_0, Z_t) \geq \frac{\varepsilon_n}{4}\right\}.$$

For every $z \in U_n$ with $\rho(Z_0, z) \leq \varepsilon_n/4$ we define a $T_z M$-valued stochastic process:

$$G(z, t) = -\int_0^t \int_M \left(\int_M \psi_{(z, F_s(x))} \, d\mu_0(x)\right)^{-1} \dot{J}(0_z, \delta F_s(y))(0) \, d\mu_0(y),$$

$$0 \leq t < T_n^1.$$

Clearly $G(z, t)$ is a semimartingale for each $z$ with bounded local characteristics depending smoothly on the spatial parameter $z$ (see, e.g., [27], Theorem (3.3)). So one can solve the equation

(2.14) $\qquad \delta Z'_t = G(Z'_t, \delta t), \qquad Z'_0 = Z_0.$

Its solution $Z'_t$ exists and is a semimartingale up to time

$$T'_n := T_n^1 \wedge \inf\left\{t > 0, \rho(Z'_t, Z_t) \geq \frac{\varepsilon_n}{4}\right\}.$$

We only need to show that $Z'_t = Z_t$, the exponential barycenter of $F_t(\mu_0)$, that is,

$$\int_M \dot{\gamma}(Z'_t, F_t(x))(0) \, d\mu_0(x) = 0_{Z'_t}. \tag{2.15}$$

The equality holds for $t = 0$. The process $\dot{\gamma}(Z'_t, F_t(x))(0)$ is a semimartingale with continuous spatial parameter and local characteristics uniformly bounded on $(x,t) \in \mathrm{supp}(\mu_0) \times [0, T'_n[$ and so using [27], Theorem (3.3), we see that integration over $M$ commutes with covariant Stratonovich differentiation $D_t$ in $t$:

$$D_t \int_M \dot{\gamma}(Z'_t, F_t(x))(0) \, d\mu_0(x) = \int_M D_t \dot{\gamma}(Z'_t, F_t(x))(0) \, d\mu_0(x). \tag{2.16}$$

More precisely, (2.16) is equivalent to: writing $V_t(x) = \dot{\gamma}(Z'_t, F_t(x))(0)$,

$$\int_0^t \left\langle p^*(\alpha), \delta\left(\int_M V_s(x) \, d\mu_0(x)\right)\right\rangle = \int_M \left(\int_0^t \langle p^*(\alpha), \delta V_s(x)\rangle\right) d\mu_0(x),$$

for every $0 \leq t < T'_n$ and every section $\alpha$ of $T^*M$. Here $p: TTM \to TM$ is the map induced by the connection which to an element of $TTM$ associates its vertical part.

Now since the connection is torsion free and since $\gamma$ is smooth, $D_t \frac{\partial}{\partial s} = \frac{D}{ds}\delta_t$ where $\frac{D}{ds}$ denotes covariant differentiation in $s$. Consequently, we have on $[0, T'_n[$,

$$\int_M D_t \dot{\gamma}(Z'_t, F_t(x))(0) \, d\mu_0(x)$$
$$= \int_M \frac{D}{ds} \delta_t \gamma(Z'_t, F_t(x))(0) \, d\mu_0(x)$$
$$= \int_M \dot{J}(\delta Z'_t, \delta F_t(x))(0) \, d\mu_0(x)$$
$$= \int_M \dot{J}(\delta Z'_t, 0_{F_t(x)})(0) + \dot{J}(0_{Z'_t}, \delta F_t(x))(0) \, d\mu_0(x)$$
$$= \left(\int_M \psi_{(Z'_t, F_t(x))}(\cdot) \, d\mu(x)\right)(\delta Z'_t) + \int_M \dot{J}(0_{Z'_t}, \delta F_t(x))(0) \, d\mu_0(x).$$

Plugging the expression (2.14) for $\delta Z'_t$ in the last formula yields

$$D_t \int_M \dot{\gamma}(Z'_t, F_t(x))(0) \, d\mu_0(x) = 0_{Z'_t}. \tag{2.17}$$

Together with $Z'_0 = Z_0$ we see that for $t \in [0, T'_n[$,

$$\int_M \dot{\gamma}(Z'_t, F_t(x))(0) \, d\mu_0(x) = 0_{Z'_t}. \tag{2.18}$$

Consequently $Z'_t$ is the exponential barycenter of $F_t(\mu_0)$ up to time $T'_n$. This implies that $T'_n = T^1_n$.

On the set $\{T^1_n < \tau_n\}$ we replace time $0$ by time $T^1_n$ and time $T^1_n$ by

$$T^2_n := \tau_n \wedge \inf\left\{t > T^1_n,\, \rho(Z_{T^1_n}, Z_t) \geq \frac{\varepsilon_n}{4}\right\},$$

to prove the requested property up to time $T^2_n$. With the same procedure we define a sequence of stopping times $(T^k_n)_{k \geq 1}$ converging almost surely to $\tau_n$, such that Proposition 2.10 is true on $[0, T^k_n[$. Consequently it is true on $[0, \tau_n[$ as requested. $\square$

2.5. *The hyperbolic space example.* Let $M$ be a Riemannian manifold and let $\nabla$ be the Levi–Civita connection. Assume that $M$ is convex for the Levi–Civita connection. Denote by $\rho(x, y)$ the distance function between $x$ and $y$ which is abbreviated as $\rho$ where there is no risk of confusion. For any two points $x \neq y$ and $\gamma = \gamma(x, y)$ the geodesic connecting them, a vector $u \in T_x M$ has an orthogonal decomposition $u = u^L + u^N$ where $u^L$ is in the direction of $\dot\gamma(0)$. We may also use $u^{L(x,y)}$ and $u^{N(x,y)}$ when the points concerned need to be clarified. Denote by $//_{x,y}: T_x M \to T_y M$ the parallel transport along the geodesic $\gamma(x, y)$. Set $E_i(x, y)(s) = (//_{x, \gamma(s)}) E_i(x, y)(0)$ for some tangent vectors $E_i(x, y)(0)$ on $T_x M$ such that $(E_1(x, y)(0), \ldots, E_n(x, y)(0))$ is an orthonormal frame of $T_x M$. We shall fix $E_1(x, y)(0) = \frac{\dot\gamma(0)}{\|\dot\gamma(0)\|}$ so that

$$E_1(x, y)(s) = \frac{\exp^{-1}_{\gamma(x,y)(s)} y}{\rho(\gamma(x, y)(s), y)}, \qquad s \in [0, 1[.$$

Since $M$ is convex we can assume that the $E_i(x, y)$'s are chosen to be smooth in $x$ and $y$ in a neighborhood of some $(x_0, y_0)$ off diagonal.

Denote by $\mathbb{H}^m$ the hyperbolic space of dimension $m$ and by $G: \mathbb{H}^m \times \mathbb{H}^m \to \mathbb{H}^m$ the smooth map which sends $(x, y)$ to the center of the geodesic segment $[x, y]: G(x, y) = \exp_x(\frac{1}{2} \exp^{-1}_x y)$.

LEMMA 2.13. *The map $G: \mathbb{H}^m \times \mathbb{H}^m \to \mathbb{H}^m$ is harmonic.*

PROOF. We need to prove that $\operatorname{tr} \nabla dG = 0$. We do it with a symmetry argument. Observe that $G(x, y) = G(y, x)$ for all $x, y$. We have

(2.19) $$\operatorname{tr} \nabla dG(x, y) = \operatorname{tr} \nabla dG(y, x).$$

On the other hand, the symmetry $\varphi_{x,y}$ in $\mathbb{H}^m$ with center $G(x, y)$ is an isometry which exchanges $x$ and $y$. So

(2.20) $$\operatorname{tr} \nabla dG(y, x) = \operatorname{tr} \nabla dG(\varphi_{x,y}(x), \varphi_{x,y}(y))$$
$$= (\varphi_{x,y})_* \operatorname{tr} \nabla dG(x, y) = -\operatorname{tr} \nabla dG(x, y).$$

The identities (2.19) and (2.20) lead to $\operatorname{tr} \nabla dG(y,x) = 0$. □

PROPOSITION 2.14. *For $x_0, y_0 \in \mathbb{H}^m$ set $\mu_0 = \frac{1}{2}(\delta_{x_0} + \delta_{y_0})$. Consider two independent Brownian motions $X_t$ and $Y_t$ with initial points $x_0$ and $y_0$, respectively. Set $\rho_t = \rho(X_t, Y_t)$ and define a measure $\mu_t$ by $\mu_t(A) = \frac{1}{2}(\delta_{X_t}(A) + \delta_{Y_t}(A))$. Then $Z_t$ is a martingale in $M$. It satisfies the Itô equation*

$$d^\nabla Z_t = \frac{1}{2}(//_{X_t,Z_t}\, d^\nabla X_t)^L + \frac{1}{2}(//_{Y_t,Z_t}\, d^\nabla Y_t)^L$$
$$+ \frac{1}{2\cosh(\rho_t/2)}((//_{X_t,Z_t}\, d^\nabla X_t)^N + (//_{Y_t,Z_t}\, d^\nabla Y_t)^N),$$

*where $u = u^L + u^N$ with $u^L$ tangential (resp. $u^N$ orthogonal) to the geodesic $\gamma(X_t, Y_t)$. Furthermore, we have locally*

$$d^\nabla Z_t = \frac{1}{\sqrt{2}} E_1(Z_t, X_t)\, dW_t^1 + \frac{1}{\sqrt{2}\cosh(\rho_t/2)} \sum_{i=2}^m E_i(Z_t, X_t)\, dW_t^i$$

*for $W_t^i$ independent real-valued Brownian motions.*

PROOF. We have $Z_t = G(X_t, Y_t)$, where $G$ is smooth and harmonic. Consequently, since $(X, Y)$ is a Brownian motion in $\mathbb{H}^m \times \mathbb{H}^m$, $Z_t$ is a martingale in $\mathbb{H}^m$. We first look for a Stratonovich equation for $Z_t$. All Jacobi fields along $\gamma(x, y)$ can be written as: $J(s) = (a^1 + b^1 s) E_1(s) + \sum_{i=2}^m [\cosh(s\rho)a^i + \sinh(s\rho)b^i] E_i(s)$, where $a_i, b_i \in \mathbb{R}$ and $\rho = \rho(x,y)$. In particular for $u \in T_x \mathbb{H}^m$, $v \in T_y \mathbb{H}^m$ and $J(u,v)$ the Jacobi fields along the geodesic $\gamma(x,y)$ with boundary values $u$ and $v$, one has

$$J(u,v)(s) = (u^1(1-s) + sv^1) E_1(s)$$
$$+ \sum_{i=2}^m \left((\cosh(s\rho) - \sinh(s\rho)\coth\rho) u^i + \frac{\sinh(s\rho)}{\sinh\rho} v^i\right) E_i(s),$$

where $u^i = \langle u, E_i(0)\rangle_\mathbb{H}$, $v^i = \langle v, E_i(1)\rangle_\mathbb{H}$ and $\langle \cdot, \cdot \rangle_\mathbb{H}$ denotes the scalar product in $\mathbb{H} = \mathbb{H}^m$. Note that $u^1 E_1(0) = u^{L(x,y)}$ and $v^1 E_1(1) = v^{L(x,y)}$. This gives

$$\dot J(u,0)(0) = -u^1 E_1(0) - \sum_{i=2}^m \rho(x,y)\coth(\rho(x,y)) u^i E_i(0)$$

(2.21)
$$= -u^{L(x,y)} - \rho(x,y)\coth(\rho(x,y)) u^{N(x,y)},$$

$$\dot J(0,v)(0) = v^1 E_1(0) + \sum_{i=2}^m \frac{\rho(x,y)}{\sinh(\rho(x,y))} v^i E_i(0)$$

$$= (//_{y,x} v)^{L(x,y)} + \frac{\rho(x,y)}{\sinh(\rho(x,y))} (//_{y,x} v)^{N(x,y)}.$$

We have

$$\int_{\mathbb{H}^m} \dot{J}(\delta Z_t, 0_{F_t(y)})(0) \, d\mu_0(y)$$
$$= -\tfrac{1}{2}(\delta Z_t)^{L(Z_t, X_t)} - \tfrac{1}{2}\rho(Z_t, X_t) \coth(\rho(Z_t, X_t))(\delta Z_t)^{N(Z_t, X_t)}$$
$$\quad - \tfrac{1}{2}(\delta Z_t)^{L(Z_t, Y_t)} - \tfrac{1}{2}\rho(Z_t, Y_t) \coth(\rho(Z_t, Y_t))(\delta Z_t)^{N(Z_t, Y_t)}$$
$$= -(\delta Z_t)^{L(Z_t, X_t)} - \tfrac{1}{2}\rho(X_t, Y_t) \coth(\tfrac{1}{2}\rho(X_t, Y_t))(\delta Z_t)^{N(Z_t, X_t)},$$

since $X_t$, $Y_t$ and $Z_t$ lie on the same geodesic $v^{L(Z_t, X_t)} = v^{L(Z_t, Y_t)}$ for any vector $v$. We shall abbreviate $v^{L(Z_t, X_t)} = v^{L(Z_t, Y_t)}$ as $v^L$. On the other hand

$$\int_{\mathbb{H}^m} \dot{J}(0_{Z_t}, \delta F_t(x))(0) \, d\mu_0(x)$$
$$= \tfrac{1}{2}\dot{J}(0, \delta X_t)(0) + \tfrac{1}{2}\dot{J}(0, \delta Y_t)(0)$$
$$= \tfrac{1}{2}(//_{X_t, Z_t}\delta X_t)^L + \tfrac{1}{2}\frac{\rho(Z_t, X_t)}{\sinh(\rho(Z_t, X_t))}(//_{X_t, Z_t}\delta X_t)^N$$
$$\quad + \tfrac{1}{2}(//_{Y_t, Z_t}\delta Y_t)^L + \tfrac{1}{2}\frac{\rho(Z_t, Y_t)}{\sinh(\rho(Z_t, Y_t))}(//_{Y_t, Z_t}\delta Y_t)^N$$
$$= \tfrac{1}{2}(//_{X_t, Z_t}\delta X_t)^L + \tfrac{1}{4}\frac{\rho(X_t, Y_t)}{\sinh(\tfrac{1}{2}\rho(X_t, Y_t))}(//_{X_t, Z_t}\delta X_t)^N$$
$$\quad + \tfrac{1}{2}(//_{Y_t, Z_t}\delta Y_t)^L + \tfrac{1}{4}\frac{\rho(X_t, Y_t)}{\sinh(\tfrac{1}{2}\rho(X_t, Y_t))}(//_{Y_t, Z_t}\delta Y_t)^N.$$

Finally, applying Proposition 2.10(ii) gives the Stratonovich equation

$$\delta Z_t = \tfrac{1}{2}(//_{X_t, Z_t}\delta X_t)^L + \tfrac{1}{2}(//_{Y_t, Z_t}\delta Y_t)^L$$
$$\quad + \tfrac{1}{2}\cosh(\tfrac{1}{2}\rho(X_t, Y_t))((//_{X_t, Z_t}\delta X_t)^N + (//_{Y_t, Z_t}\delta Y_t)^N).$$

Obtaining the Itô equation (2.14) from the Stratonovich equation is immediate since $Z_t$ is a martingale.

For the local form of the Itô equation, by a localization procedure as in [8], Lemma 3.5, we only need to define $W^i$ locally. We do it with the formula

$$(2.22) \qquad dW_t^i = \frac{1}{\sqrt{2}}\langle E_i(Z_t, X_t), //_{X_t, Z_t} d^\nabla X_t + //_{Y_t, Z_t} d^\nabla Y_t\rangle_{\mathbb{H}},$$

valid for $(X_s, Y_s)$ in a small open subset of $\mathbb{H}^m \times \mathbb{H}^m \setminus \{\text{diagonal}\}$. The processes $(W_t^i)$ are martingales with quadratic variation $t$ and therefore Brownian motions. The independence comes from the mutual orthogonality of $E_i(Z_t, X_t)$'s. □

REMARK 2.15. Note that if $\omega_0$ is a point such that the geodesic $(X_t(\omega_0)Y_t(\omega_0))$ converges, then $Z_t$ is close to the solution to $dz_t = \frac{1}{\sqrt{2}}E_1(z_t, X_\infty(\omega_0))\, dW_t^1$, which lives on the geodesic segment $(X.(\omega_0)Y.(\omega_0))_\infty$ if it starts there. In the next section we shall make precise this convergence.

**3. The barycenter of two independent Brownian particles in $\mathbb{H}$.** Consider the hyperbolic space $\mathbb{H} = \mathbb{H}^d$ of dimension $d \geq 2$ as an oriented manifold. Let $\partial \mathbb{H}$ be the visibility boundary in its visibility compactification, that is, the set of equivalence classes of geodesic rays. In the upper half space model $\partial \mathbb{H}$ is the union of the boundary hyperplane and the point at infinity. The visibility topology coincides with the usual topology when $\mathbb{H}$ is considered as a subset of $\mathbb{R}^d$.

Let $X$ and $Y$ be two independent Brownian motions in $\mathbb{H}$ with initial points $X_0$ and $Y_0$, respectively. Then they converge in $\bar{\mathbb{H}}$ with limits $X_\infty$ and $Y_\infty$ in $\partial \mathbb{H}$, respectively. We denote by $(XY)_\infty \subset \mathbb{H}$ the random geodesic connecting the two limit points. By a random variable $u$ on $(XY)_\infty$ we mean $u(\omega) \in (XY)_\infty(\omega)$ for almost all $\omega$.

THEOREM 3.1. *Let $X_t$ and $Y_t$ be two independent Brownian particles in $\mathbb{H}$ and let $Z_t$ be their exponential barycenter. Denote by $\mu(X_0, Y_0)$ the law of $(Z_t)_{t \geq 0}$. Let $u$ be an $\mathcal{F}_\infty$-measurable random variable on $(XY)_\infty$ and let $(T_n)_{n \geq 1}$ be a sequence of finite stopping times increasing to infinity such that $Z_{T_n}$ converges to $u$ almost surely. Then for almost all $\omega_0$, $\mu(X_{T_n}(\omega_0), Y_{T_n}(\omega_0))$ converges as $n$ goes to infinity to the law of a Brownian motion $(z_t)_{t \geq 0}$ of variance $1/2$ on the geodesic $(XY)_\infty(\omega_0)$ with starting point $u(\omega_0)$; that is, $(z_{2t})_{t \geq 0}$ is a Brownian motion on $(XY)_\infty(\omega_0)$.*

PROOF. For the proof take the upper half space model
$$\mathbb{H} \equiv \{y = (y^1, y^2, \ldots, y^d) \in \mathbb{R}^d, y^d > 0\}.$$
The process $(X_t, Z_t)$ is a diffusion as the image of the diffusion $(X, Y)$ by the diffeomorphism $(x, y) \mapsto (x, \gamma(x,y)(\frac{1}{2}))$ on $\mathbb{H} \times \mathbb{H}$. Set $X_t^n = X_{T_n+t}$, $Y_t^n = Y_{T_n+t}$ and $Z_t^n = Z_{T_n+t}$. Then $(X_t^n, Z_t^n)$ has initial condition $(X_{T_n}, Z_{T_n})$ which by assumption converges to $(X_\infty, u)$.

Since the hyperbolic Laplacian is given by $\Delta_\mathbb{H} f(x) = \frac{1}{2}(x^d)^2 \Delta - \frac{d-2}{2}x^d \partial_d$, the hyperbolic Brownian motions can be written as the solution to the Itô stochastic differential equation

(3.1)
$$dX_t = X_t^d\, dB_t - \tfrac{1}{2}(d-2)X_t^d e_d\, dt,$$
$$dY_t = Y_t^d\, dB_t' - \tfrac{1}{2}(d-2)Y_t^d e_d\, dt,$$

where $(B_t, B_t')$ is a Brownian motion in $\mathbb{R}^{2d}$ and $\{e_i\}_{i=1}^d$ is the canonical basis in $\mathbb{R}^d$.

Set $f(x^1, \ldots, x^d) = \sup(x^d, 0)$. Then $(X_t)$ can be considered to be the solution to
$$dX_t = f(X_t)\, dB_t - \tfrac{1}{2}(d-2)f(X_t)e_d\, dt$$
on $\mathbb{R}^d$. This equation has locally Lipschitz coefficients. We shall show that the process $(X_t, Z_t)$ also satisfies an equation which extends to $\mathbb{R}^d \times \mathbb{H}$ with local Lipschitz coefficients. First note that since $X_t$ and $Y_t$ are martingales, (3.1) yields
$$d^\nabla X_t = X_t^d\, dB_t, \qquad d^\nabla Y_t = Y_t^d\, dB'_t.$$
On the other hand, by Proposition 2.14
$$d^\nabla Z_t = \frac{1}{2}(/\!/_{X_t, Z_t} d^\nabla X_t)^L + \frac{1}{2}(/\!/_{Y_t, Z_t} d^\nabla Y_t)^L$$
$$+ \frac{1}{2\cosh(\rho_t/2)}((/\!/_{X_t, Z_t} d^\nabla X_t)^N + (/\!/_{Y_t, Z_t} d^\nabla Y_t)^N).$$
We have
$$(3.2) \qquad dZ_t = d^\nabla Z_t - \tfrac{1}{2}\Gamma(Z_t)(d^\nabla Z_t, d^\nabla Z_t)\, dt,$$
where $\Gamma(z)$ is the Christoffel symbol at $z$ for the canonical connection in $\mathbb{H}$. It is a smooth function, so we only need to prove that the equation for $d^\nabla Z_t$ extends as requested. For this plug $d^\nabla X_t$ and $d^\nabla Y_t$ into the formula for $d^\nabla Z_t$ and observe that $Y(t) = H(X(t), Z(t))$ for $H(x, z) = \exp(2\exp_x^{-1} z)$ to obtain
$$d^\nabla Z_t = \frac{1}{2}(/\!/_{X_t, Z_t} (X_t^d dB_t)^{L(X_t, Z_t)}$$
$$+ /\!/_{H(X_t, Z_t), Z_t} (H(X_t, Z_t)^d\, dB'_t)^{L(H(X_t, Z_t), Z_t)})$$
$$+ \frac{/\!/_{X_t, Z_t}(X_t^d dB_t)^{N(X_t, Z_t)}}{2\cosh(\rho(X_t, Z_t))}$$
$$+ \frac{/\!/_{H(X_t, Z_t), Z_t}(H(X_t, Z_t)^d\, dB'_t)^{N(H(X_t, Z_t), Z_t)}}{2\cosh(\rho(X_t, Z_t))}.$$
Let $w = (w^1, \ldots, w^d)$ be a vector in $\mathbb{R}^d$ and let $x_n = (x_n^1, \ldots, x_n^d)$ be a sequence in $\mathbb{H}$ converging to $x = (x^1, \ldots, x^{d-1}, 0) \in \partial\mathbb{H} \setminus \{\infty\}$. Set $E(z, x) = \frac{1}{\rho(z,x)} \exp_z^{-1} x$, which extends smoothly to the set $(\mathbb{H} \times \overline{\mathbb{H}}) \setminus \{\text{diagonal}\}$ considered as a subset of $\mathbb{H} \times \mathbb{R}^d$. Furthermore,
$$\lim_{n\to\infty} /\!/_{x_n, z}(x_n^d w)^{L(x_n, z)}$$
$$= \lim_{n\to\infty} (-x_n^d)\langle w, E(x_n, z)\rangle_\mathbb{H} E(z, x_n)$$
$$= \lim_{n\to\infty} \langle w, -(x_n^d)^{-1} E(x_n, z)\rangle E(z, x_n) = -w^d E(z, x)$$

since $(x_n^d)^{-1}E(x_n,z)$ converges to $e_d$, and similarly

$$\lim_{n\to\infty} //_{H(x_n,z),z}(H(x_n,z)^d w)^{L(H(x_n,z),z)} = w^d E(z,x).$$

Moreover, for $z \in \mathbb{H}$,

$$\lim_{n\to\infty} \frac{(//_{x_n,z}(x_n^d w)^{N(x_n,z)} + //_{H(x_n,z),z}(H(x_n,z)^d w)^{N(H(x_n,z),z)})}{2\cosh \rho(x_n,z)} = 0,$$

from the boundedness of the nominator and the convergence to $+\infty$ of $\cosh \rho(x_n,z)$. We conclude that each coefficient in the equation for $(X_t, Z_t)$ smoothly extends over $\bar{\mathbb{H}} \times \mathbb{H}$ and $(X_t, Z_t)$ is a solution to the system of equations of the following form:

(3.3) $$dX_t = f(X_t^d)\,dB_t - \frac{d-2}{2}f(X_t^d)e_d\,dt,$$

(3.4) $$dZ_t = \sigma(X_t, Z_t)\,d(B_t, B_t') + b(X_t, Z_t)\,dt,$$

where $\sigma(x,z) = \sigma(\bar{x},z)$ and $b(x,z) = b(\bar{x},z)$ if $x^d < 0$ and $\bar{x} = (x^1,\ldots,x^{d-1},0)$. Since all coefficients are smooth on $\bar{\mathbb{H}} \times \mathbb{H}$ and $(\mathbb{R}^d \setminus \mathbb{H}) \times \mathbb{H}$ they are locally Lipschitz on $\mathbb{R}^d \times \mathbb{H}$. Consequently the system (3.3)–(3.4) has a unique solution and it does not explode for starting points in $\mathbb{H} \times \mathbb{H}$ as known. If, however, the starting point satisfies $X_0^d \leq 0$, then $X_\cdot \equiv X_0$ and (3.4) reduces to

(3.5) $$d^\nabla Z_t = \tfrac{1}{2}E(Z_t, \bar{X}_0)(-dB_t^d + dB_t'^d)$$

whose solution starting from $Z_0$ shall be denoted by $(Z_t(X_0, Z_0))$. Then $Z_{2t}(X_0, Z_0)$ is a Brownian motion on the hyperbolic geodesic $(\bar{X}_0 Z_0)$. In particular, it does not explode.

Consider the process $(X^n, Z^n)$ as solution to (3.3)–(3.4) where $(B_t, B_t')$ is replaced by $(B_{T_n+t}, B'_{T_n+t})$, with starting point $(X_{T_n}, Z_{T_n})$ which by assumption converges almost surely to $(X_\infty, u)$. Since $\mathbb{R}^d \times \mathbb{H}$ is diffeomorphic to $\mathbb{R}^{2d}$ we can apply Corollary 11.1.5 in [25] for $\mathbb{R}^m$-valued SDE and conclude that the transition probabilities are Feller continuous and so by the Markov property the law of $(X^n, Z^n)$ conditioned by $(X_{T_n}, Z_{T_n})$ converges almost surely to the law of $(X_\infty, Z(X_\infty, u))$. Hence the law of $Z_{2t}$ is that of a Brownian motion on the geodesic $(X_\infty, u)$, as requested. □

The rest of this section is devoted to the existence of sequences of stopping times $(T_n)_{n\geq 1}$ as in Theorem 3.1.

First we note the following fact: If $X_t$ and $Y_t$ are independent Brownian motions in $\mathbb{H}$, then for almost all $\omega$ the random geodesic $(X_t(\omega)Y_t(\omega))$ converges to $(X(\omega)Y(\omega))_\infty$ uniformly on compact sets as $t$ goes to infinity, that is, for any compact set $K$, $\lim_{t\to\infty} \sup_{z\in(X_t Y_t)\cap K} \rho(z, (X_\infty Y_\infty)) = 0$.

To see this we take the upper half space representation. We may assume $K \subset \{(x^1, \ldots, x^d) \in \mathbb{H}, \; x^d > \varepsilon\}$ for some $\varepsilon > 0$. If $z \in (X_t Y_t) \cap K$, then $z$ belongs to a Euclidean tubular neighborhood of the half-circle $(XY)_\infty$, with radius $\sup_{s \geq t} \max(\|X_s - X_\infty\|, \|Y_s - Y_\infty\|)$. Since for $z \in K$ the hyperbolic distance is smaller than the Euclidean distance divided by $\varepsilon$,

$$\lim_{t \to \infty} \sup_{z \in (X_t Y_t) \cap K} \rho(z, (X_\infty Y_\infty)) \leq \lim_{t \to \infty} \frac{1}{\varepsilon} \sup_{s \geq t} \left( \max(\|X_s - X_\infty\|, \|Y_s - Y_\infty\|) \right)$$
$$= 0.$$

Next for $x, y, x \in \mathbb{H}$ with $x \neq y$, denote by $(xy)$ the geodesic segment connecting them and by $p(z, x, y)$ the orthogonal projection of $z$ into the geodesic $(xy)$.

LEMMA 3.2. *Let $(X_t)$, $(Y_t)$, $(T_n)_{n \geq 1}$ and $u$ be as in Theorem 3.1. If $(u_n)_{n \geq 1}$ is a sequence of $(\mathcal{F}_{T_n})$-adapted random variables in $\mathbb{H}$ converging almost surely to an $\mathcal{F}_\infty$-measurable random measurable $u$ on $(XY)_\infty$, then $\lim_{n \to \infty} p(u_n, X_{T_n}, Y_{T_n}) = u$ almost surely.*

PROOF. Since the projection to the convex set $(X_{T_n} Y_{T_n})$ is 1-Lipschitz we have

$$\rho(u, p(u_n, X_{T_n}, Y_{T_n}))$$
$$\leq \rho(u, p(u, X_{T_n}, Y_{T_n})) + \rho(p(u, X_{T_n}, Y_{T_n}), p(u_n, X_{T_n}, Y_{T_n}))$$
$$\leq \rho(u, (X_{T_n} Y_{T_n})) + \rho(u, u_n) \to 0$$

following from the fact that for almost all $\omega$ the geodesic $(X_t(\omega) Y_t(\omega))$ converges to $(XY)_\infty(\omega)$ uniformly on compact sets $K$ of $\mathbb{H}$. □

THEOREM 3.3. *Let $(X_t)$, $(Y_t)$, $(Z_t)$ and $u$ be as in Theorem 3.1. There exists a sequence of finite stopping times $(T_n)_{n \geq 1}$ increasing to infinity such that $Z_{T_n}$ converges to $u$ almost surely.*

PROOF. Let $(u_n)_{n \geq 1}$ be a sequence of $(\mathcal{F}_n)$-adapted random variables in $\mathbb{H}$ converging almost surely to $u$. By Lemma 3.2 it is sufficient to prove that there exists a sequence of stopping times $(T_n)_{n \geq 1}$ such that $T_n \geq n$ and $Z_{T_n} = p(u_n, X_{T_n}, Y_{T_n})$ almost surely. Conditioning with respect to $\mathcal{F}_n$, it is sufficient to prove that for every $o \in M$ and $n \geq 1$, there exists a stopping time $T_n \geq n$ such that $Z_{T_n} = p(o, X_{T_n}, Y_{T_n})$.

Let $P_t = p(o, X_t, Y_t)$. Since $P_t$ and $Z_t$ belong to the geodesic segment $[X_t, Y_t]$, $P_t = Z_t$ if and only if the signed distance $D_t \equiv \rho(P_t, X_t) - \rho(X_t, Z_t)$ between $P_t$ and $Z_t$ is zero. For the existence of the stopping times we only need to show that $D_t$ is recurrent. Set $R_t^X = \rho(o, X_t)$ and $R_t^Y = \rho(o, Y_t)$. By the definition of $P_t$ and the triangle inequalities, we have $\rho(X_t, P_t) \leq R_t^X$ and $\rho(P_t, Y_t) \geq R_t^Y - \rho(o, P_t)$. So

$$D_t = \rho(P_t, X_t) - \tfrac{1}{2}\rho(X_t, Y_t) = \rho(P_t, X_t) - \tfrac{1}{2}(\rho(X_t, P_t) + \rho(P_t, Y_t))$$
$$= \tfrac{1}{2}\rho(X_t, P_t) - \tfrac{1}{2}\rho(P_t, Y_t) \leq \tfrac{1}{2}(R_t^X - R_t^Y) + \tfrac{1}{2}\rho(o, P_t)$$

and similarly $D_t \geq \tfrac{1}{2}(R_t^X - R_t^Y) - \tfrac{1}{2}\rho(o, P_t)$. Since $\lim_{t\to\infty} \rho(o, P_t)$ exists we only need to show $R_t^X - R_t^Y$ is recurrent. Note that

$$dR_t^X = dB_t^X + \frac{d-1}{2}\coth R_t^X\, dt,$$

$$dR_t^Y = dB_t^Y + \frac{d-1}{2}\coth R_t^Y\, dt,$$

where $(B_t^X, B_t^Y)$ is a planar Brownian motion. So for almost all $\omega$,

$$\inf(R_t^X(\omega), R_t^Y(\omega)) > \tfrac{1}{4}(d-1)t$$

for sufficiently large time. Now

$$R_t^X - R_t^Y = B_t^X - B_t^Y + \tfrac{1}{2}(d-1)\int_0^t (\coth R_s^X - \coth R_s^Y)\, ds.$$

But $B_t^X - B_t^Y$ is recurrent and $\int_0^\infty (\coth R_s^X - \coth R_s^Y)\, ds$ exists since

$$|\coth R_s^X - \coth R_s^Y| = |\coth R_s^X - 1 - (\coth R_s^Y - 1)|$$
$$= \left|\frac{2}{e^{2R_s^X} - 1} - \frac{2}{e^{2R_s^Y} - 1}\right| \leq \frac{4}{e^{(d-1)s/2} - 1}$$

for large time $s$. Consequently $R_t^X - R_t^Y$ is recurrent and so is $D_t$. □

**4. Barycenters of measures in a Cartan–Hadamard manifold.** Let $M$ be a Cartan–Hadamard manifold with pinched negative curvatures $k$, $-b^2 \leq k \leq -a^2$, for $b \geq a > 0$. Denote by $\bar M = M \cup \partial M$ its visibility compactification where $\partial M$ is the set of equivalence classes of unit speed geodesics in $M$ under the equivalence relation

$$\gamma_1 \sim \gamma_2 \iff \limsup_{t\to\infty} \rho(\gamma_1(t), \gamma_2(t)) < \infty,$$

endowed with the sphere topology. See, for example, [2], page 22. Note that for each point $z \in M$ and $x \in \partial M$ there is a unique unit speed geodesic in the equivalence class of $x$ with initial point $z$, which we shall denote as $\{\varphi(z, x)(t),\ 0 \leq t < \infty\}$. For $z, y \in M$, $z \neq y$, we shall also denote by $\{\varphi(z, y)(t), t \geq 0\}$ the unit speed geodesic satisfying $\varphi(z, y)(0) = z$ and $\varphi(z, y)(\rho(z, y)) = y$. In other words,

$$\varphi(z, y)(t) = \gamma(z, y)(t/\rho(z, y)) \qquad \forall z, y \in M \text{ satisfying } z \neq y.$$

Let $\pi: TM \to M$ be the canonical projection and let

$$\phi: SM \equiv \{v \in TM, \ \|v\| = 1\} \longrightarrow \partial M$$

be the map which sends $\theta$ to $[\exp_{\pi(\theta)}(t\theta)]_{t\geq 0} \in \partial M$. The map $\Phi = (\pi, \phi)$ from $SM$ to $M \times \partial M$ is a homeomorphism. In fact, $\Phi^{-1}(z,x) = \dot\varphi(z,x)(0)$ (see, e.g., [6], Propositions 2.13 and 2.14).

For $o \in M$ denote by $(\psi_{o,x}, x \in \partial M)$ the family of Busemann functions:

(4.1) $$\psi_{o,x}(y) = \lim_{t \to \infty} [\rho(y, \varphi(o,x)(t)) - t].$$

These functions are characterized by $\psi_{o,x}(o) = 0$ and $\mathrm{grad}_y \psi_{o,x} = -\Phi^{-1}(y,x) = -\dot\varphi(y,x)(0)$, any $y \in M$. We write $\psi_x$ for $\psi_{o,x}$ if there is no risk of confusion. Denote by $\mathcal{M}_1(X)$ the set of probability measures on a topological space $X$ endowed with the Borel $\sigma$-field. Set

(4.2) $$\psi_\mu(z) \equiv \psi_{o,\mu}(z) = \int_{\partial M} \psi_{o,x}(z) \, d\mu(x), \qquad \mu \in \mathcal{M}_1(\partial M).$$

Then

(4.3) $$\mathrm{grad}_z \psi_\mu = -\int_{\partial M} \dot\varphi(z,x)(0) \, d\mu(x).$$

A solution to $\mathrm{grad}_z \psi_\mu = 0$ is called a Busemann barycenter of $\mu$.

The aim of this section is to show that for discrete probability measures of finite support transported by a suitable random flow, the exponential barycenter $\mu_t$ converges to the Busemann barycenter of the limiting measure on the boundary.

LEMMA 4.1. *Let $\mu$ be a probability measure on $M$ with corresponding vector field $H_\mu : M \to TM$ defined by (2.1). We have*

$$\langle \nabla_u H_\mu, u \rangle \leq -a\varepsilon(z,\mu) \|u\|^2, \qquad u \in T_z \mathbb{H},$$

*where*

(4.4) $$\varepsilon(z,\mu) = \min_{w \in S_z M} \int_M \rho(z,x) \sin^2(w, \dot\varphi(z,x)(0)) \, d\mu(x).$$

PROOF. We follow the notation of Section 2.3 and note that

$$\nabla_u H_\mu = \int_M \dot J(u, 0_x)(0) \, d\mu(x),$$

where $0_x \in T_x M$ is the zero vector. Write $u = u^L + u^N$ where $u^L \equiv u^{L(z,x)}$ is colinear to $\exp_z^{-1} x$ and $u^N \equiv u^{N(z,x)}$ is its orthogonal complement. Then

$$\dot J(u, 0_x)(0) = \dot J(u^L, 0_x)(0) + \dot J(u^N, 0_x)(0) = -u^L + \dot J(u^N, 0_x)(0),$$

sum of two orthogonal vectors. Write $J(t) = J(u^N, 0_x)(t)$ and $f(t) = \|J(t)\|$. Then $f'(t) = \frac{1}{\|J(t)\|}\langle J(t), \dot{J}(t)\rangle$. In particular, $f'(0) = \frac{1}{\|u^N\|}\langle u^N, \dot{J}(0)\rangle$ and furthermore

$$f''(t) = f^{-3}(t)(\|J(t)\|^2\|\dot{J}(t)\|^2 - \langle J(t), \dot{J}(t)\rangle^2 + \|J(t)\|^2\langle J(t), \ddot{J}(t)\rangle)$$
$$\geq f^{-1}(t)\langle J(t), \ddot{J}(t)\rangle$$
$$= -f^{-1}(t)\langle J(t), R(J(t), \dot{\gamma}(t))\dot{\gamma}(t)\rangle \geq a^2\rho^2(z,x)f(t),$$

where $\gamma(t) = \gamma(z,x)(t)$ for $t \in [0,1]$ and $a$ is the upper bound of the sectional curvature. Note that $f(0) = \|u^N\|$ and $f(1) = 0$. By comparison with the solution to $g''(t) = a^2\rho^2(z,x)g(t)$ with the same boundary conditions, we see

$$f'(0) \leq g'(0) = -a\rho(z,x)\coth(a\rho(z,x))\|u^N\| \leq -a\rho(z,x)\|u^N\|$$

which leads to $\langle u^N, \dot{J}(0)\rangle \leq -a\rho(z,x)\|u^N\|^2$. Finally

$$\langle \nabla_u H_\mu, u\rangle = \int_M (-\|u^L\|^2 + \langle \dot{J}(0), u^N\rangle)\,d\mu(x)$$
$$\leq \int_M (-\|u^L\|^2 - a\rho(z,x)\|u^N\|^2)\,d\mu(x) \leq -\int_M a\rho(z,x)\|u^N\|^2\,d\mu(x)$$
$$= -a\|u\|^2 \int_M \rho(z,x)\sin^2(u, \dot{\varphi}(z,x)(0))\,\mu(dx) \leq -a\|u\|^2\varepsilon(z,\mu). \quad\square$$

LEMMA 4.2. *Let $\mu$ be a probability measure on $\partial M$. Then*

(4.5) $$\nabla d\psi_\mu(u,u) \geq a\alpha(z,\mu)\|u\|^2 \qquad \forall u \in T_zM,\ z \in M,$$

*where*

(4.6) $$\alpha(z,\mu) = \min_{w \in S_zM} \int_{\partial M} \sin^2(w, \dot{\varphi}(z,x)(0))\,d\mu(x).$$

PROOF. Let $z_0 \in M$ and $x \in \partial M$. Set

$$\psi_{t,x}(-) = \rho(-, \varphi(z_0,x)(t)) - t.$$

By Proposition 3.1 in [13], as $t$ goes to infinity, $(\psi_{t,x}, \operatorname{grad}\psi_{t,x}, \nabla\operatorname{grad}\psi_{t,x})$ converges to $(\psi_x, \operatorname{grad}\psi_x, \nabla\operatorname{grad}\psi_x)$ uniformly on compact sets. In fact, the proof of the same proposition shows that the convergence is uniform in $x$. In particular, if we set

$$\psi_{t,\mu}(z) = \int_{\partial M} \psi_{t,x}(z)\,d\mu(x),$$

then

$$\operatorname{grad}_z \psi_{t,\mu} = -\int_{\partial M} \dot{\varphi}(z, \varphi(z_0,x)(t))(0)\,d\mu(x),$$

$$\nabla_u \operatorname{grad}\psi_{t,\mu} = -\int_{\partial M} \nabla_u\dot{\varphi}(\cdot, \varphi(z_0,x)(t))(0)\,d\mu(x),$$

$(\psi_{t,\mu}, \operatorname{grad} \psi_{t,\mu}, \nabla \operatorname{grad} \psi_{t,\mu})$ converges to $(\psi_\mu, \operatorname{grad} \psi_\mu, \nabla \operatorname{grad} \psi_\mu)$ uniformly on compact sets and the convergence is uniform in $\mu$.

Let $\mu_t = \varphi(z_0, -)(t)_*(\mu)$. Since $\dot\gamma(x,y)(0) = \rho(x,y)(o)\dot\varphi(x,y)(0)$,

$$H_{\mu_t}(z) = \int_M \dot\gamma(z,y)(0) \, d\mu_t(y)$$
$$= \int_{\partial M} \rho(z, \varphi(z_0,x)(t)) \dot\varphi(z, \varphi(z_0,x)(t))(0) \, d\mu(x).$$

If $u \in T_{z_0} M$, then $d_u \rho(\cdot, \varphi(z_0,x)(t)) = -\langle \dot\varphi(z_0,x)(0), u\rangle$ and

$$\nabla_u H_{\mu_t} = \int_{\partial M} \rho(z_0, \varphi(z_0,x)(t)) \nabla_u \dot\varphi(\cdot, \varphi(z_0,x)(t))(0) \, d\mu(x)$$
$$+ \int_{\partial M} d_u \rho(\cdot, \varphi(z_0,x)(t)) \dot\varphi(z_0, \varphi(z_0,x)(t))(0) \, d\mu(x)$$
$$= \int_{\partial M} t \nabla_u \dot\varphi(\cdot, \varphi(z_0,x)(t))(0) \, d\mu(x)$$
$$- \int_{\partial M} \langle \dot\varphi(z_0,x)(0), u\rangle \dot\varphi(z_0,x)(0) \, d\mu(x)$$
$$= -t \nabla_u \operatorname{grad} \psi_{t,\mu} - \int_{\partial M} \langle \dot\varphi(z_0,x)(0), u\rangle \dot\varphi(z_0,x)(0) \, d\mu(x)$$

using $\rho(z_0, \varphi(z_0,x)(t)) = t$ and for all $s \geq 0$, $\varphi(z_0, \varphi(z_0,x)(t))(s) = \varphi(z_0,x)(s)$. Consequently

$$\langle \nabla_u H_{\mu_t}, u\rangle = -t\langle \nabla_u \operatorname{grad} \psi_{t,\mu}, u\rangle - \int_{\partial M} \langle \dot\varphi(z_0,x)(0), u\rangle^2 \, d\mu(x).$$

This together with Lemma 4.1 gives

$\langle \nabla_u \operatorname{grad} \psi_{t,\mu}, u\rangle$
$$\geq \frac{a}{t} \varepsilon(z_0, \mu_t) \|u\|^2 - \frac{1}{t} \int_{\partial M} \langle \dot\varphi(z_0,x)(0), u\rangle^2 \, d\mu(x)$$
$$\geq \frac{a}{t} \|u\|^2 \min_{w \in S_{z_0} M} \int_M \rho(z_0, y) \sin^2(w, \dot\varphi(z_0,y)(0)) \, d\mu_t(y) - \frac{1}{t}\|u\|^2$$
$$\geq \frac{a}{t} \|u\|^2 \min_{w \in S_{z_0} M} \int_{\partial M} \rho(z_0, \varphi(z_0,x)(t))$$
$$\times \sin^2(w, \dot\varphi(z_0, \varphi(z_0,x)(t))(0)) \, d\mu(x) - \frac{1}{t}\|u\|^2$$
$$= a\|u\|^2 \min_{w \in S_{z_0} M} \int_{\partial M} \sin^2(w, \dot\varphi(z_0,x)(0)) \, d\mu(x) - \frac{1}{t}\|u\|^2$$
$$= a\alpha(z_0, \mu)\|u\|^2 - \frac{1}{t}\|u\|^2.$$

Taking $t$ to infinity and using the convergence of $\nabla \operatorname{grad} \psi_{t,\mu}$ to $\nabla \operatorname{grad} \psi_\mu$, we obtain $\nabla d\psi_\mu(u,u) \geq a\alpha(z_0,\mu)\|u\|^2$, as desired. $\square$

Let $\mathcal{U}$ be the subset of $\mathcal{M}_1(\partial M)$ containing discrete measures with no atoms of weight greater than or equal to $1/2$. Clearly the function $\alpha(z,\mu)$ is continuous on $M \times \mathcal{U}$, using the weak convergence topology on $\mathcal{U}$. Moreover, it is strictly positive on $M \times \mathcal{U}$: for every $z \in M$, the set $\{\dot\varphi(z,x)(0),\, x \in \operatorname{supp}(\mu)\}$ contains at least three different vectors, and this implies that no $w \in S_z M$ is colinear to all of them. The positivity then follows from the compactness of $S_z M$.

PROPOSITION 4.3. *For $\mu \in \mathcal{U}$ there exists a unique $z \in M$ such that*

$$(4.7) \qquad \operatorname{grad}_z \psi_\mu \equiv -\int_{\partial M} \dot\varphi(z,x)(0)\,d\mu(x) = 0.$$

*Denote the solution by $G(\mu)$. Then the map $G:\mathcal{U} \to M$ is continuous.*

PROOF. The existence and uniqueness are well known in the case that $\mu$ is a continuous measure. In fact the uniqueness follows from Lemma 4.2 since $\alpha(z,\mu) > 0$ and $\psi_\mu$ is strictly convex. For the existence we only need to show that there is a geodesic ball $B(o,T) \subset M$ of radius $T > 0$ on which $\operatorname{grad} \psi_\mu$ points outward the boundary and therefore has a zero inside the geodesic ball. For $T > 0$, take $y \in \partial B(o,T)$. Let $\bar\gamma \in \partial M$ be the point corresponding to the geodesic ray $\gamma(o,y)$ and let $B_\varepsilon(\bar\gamma) \subset \partial M$ be the set of points whose angle with $\gamma(o,y)$ is smaller than $\varepsilon$ using the sphere topology. Choose $\varepsilon_0 > 0$ so that $c_0 \equiv \sup_{x \in \partial M} \mu(B_{\varepsilon_0}(x)) < 1/2$ which is possible due to the compactness of the sphere at infinity. Let $x \in \partial M$, $\varepsilon$ the angle between $\dot\varphi(0,y)(0)$ and $\dot\varphi(0,x)(0)$, and $\alpha$ the angle between $\dot\varphi(y,x)(0)$ and $\dot\varphi(y,o)(0)$. By comparison with a manifold with constant curvature $-a^2$, we have

$$\sin\alpha \leq \frac{\cos\varepsilon\cos\alpha + 1}{\sin\varepsilon \cosh(aT)} \leq \frac{2}{\sin\varepsilon \cosh(aT)}$$

and

$$\langle \dot\varphi(y,x)(0), \dot\varphi(y,o)(0)\rangle \geq 1 - \frac{4}{\sin^2\varepsilon \cosh^2(aT)}.$$

Consequently

$$-\langle \operatorname{grad}_y \psi_\mu, \dot\varphi(y,o)(0)\rangle$$
$$= \left(\int_{\partial M \setminus B_{\varepsilon_0}(\bar\gamma)} + \int_{B_{\varepsilon_0}(\bar\gamma)}\right) \langle \dot\varphi(y,x)(0), \dot\varphi(y,o)(0)\rangle\, d\mu(x)$$
$$\geq \int_{\partial M \setminus B_{\varepsilon_0}(\bar\gamma)} \left(1 - \frac{4}{\sin^2\varepsilon \cosh^2(aT)}\right) d\mu(x) + (-1)\cdot\mu(B_{\varepsilon_0}(\bar\gamma))$$

$$= \left(1 - \frac{4}{\sin^2 \varepsilon \cosh^2(aT)}\right)(1 - \mu(B_{\varepsilon_0}(\bar{\gamma}))) - \mu(B_{\varepsilon_0}(\bar{\gamma}))$$

$$\geq 1 - \frac{4}{\sin^2 \varepsilon \cosh^2(aT)} - 2c_0$$

$$> 0$$

when $\cosh^2(aT) > 4/\sin^2 \varepsilon (1 - 2c_0)$. Since $\dot{\varphi}(y,o)(0)$ is normal to the boundary $\partial B(o,T)$ and pointing inwards, this proves that $\operatorname{grad}_y \psi_\mu$ is pointing outwards. Consequently there exists $z \in B(o,T)$ such that $\operatorname{grad}_z \psi_\mu = 0$, so $G$ is well defined.

For the continuity of $G$ note that $\operatorname{grad}_y \psi_\mu$ is continuous with respect to $(y, \mu)$ since $\operatorname{grad}_y \psi_{t,\mu}$, which is continuous in $(y, t, \mu)$, converges to $\operatorname{grad}_y \psi$ uniformly on $y$ in compact sets and uniformly in $\mu$. Let $(\mu_n)_{n\geq 1}$ be a sequence of elements of $\mathcal{U}$ converging to $\mu \in \mathcal{U}$. Set $z_n = G(\mu_n)$ and $z = G(\mu)$. From the convergence of the sequence $(\mu_n)_{n\geq 1}$ to $\mu$, we can choose the same $\varepsilon_0$, $c_0$ and $T$ for all $\mu_n$ and $\mu$. Consequently, $z$ and all the $z_n$ belong to $B(o,T)$. Furthermore, the function $\alpha$ in Lemma 4.2 is continuous, so we have $\alpha_0 = \inf\{\alpha(y, \mu_n) : (y, n) \in B(o,T) \times \mathbb{N}\} > 0$ and consequently $\langle \nabla_u \operatorname{grad} \psi_{\mu_n}, u \rangle \geq a\alpha_0 \|u\|^2$ for all $y \in B(o,T)$ and $u \in T_y M$. Let $z_n(t)_{t\geq 0}$ be $C^1$ paths in $M$ with $z_n(0) = z$ and $\dot{z}_n(t) = -\operatorname{grad} \psi_{\mu_n}(z_n(t))$. Then by differentiating $\|\operatorname{grad} \psi_{\mu_n}(z_n(t))\|^2$ with respect to $t$ we see $\|\operatorname{grad} \psi_{\mu_n}(z_n(t))\| \leq e^{-a\alpha_0 t}\|\operatorname{grad}\psi_{\mu_n}(z_n(0))\|$. It follows that $z_n = z_n(\infty) = \lim_{t\to\infty} z_n(t)$ and the length of the curve from $z = z_n(0)$ to $z_n$ is smaller than $\|\operatorname{grad}\psi_{\mu_n}(z)\|/(a\alpha_0)$. Since $\operatorname{grad}\psi_{\mu_n}(z) \to \operatorname{grad}\psi_\mu(z) = 0$, $\lim_{n\to\infty} \rho(z, z_n) = 0$. So $G$ is continuous. □

EXAMPLE 4.4. Let $\mathbb{H} \subset \mathbb{C}$ be the Poincaré upper half plane. The boundary of $\mathbb{H}$ is the real line $\mathbb{R}$, plus one point at infinity. For $n \geq 2$ and $x_1 < \cdots < x_n$ in $\mathbb{R}$ set $\mu = \frac{1}{n}\sum_{j=1}^n \delta_{x_j}$. If $n = 3$, the angles between the vectors $\dot{\varphi}(G(\mu), x_j)(0)$ are $\pm 2\pi/3$. Considering first the case $x_3 = \infty$ and then the general case via homographic transformations, one finds

$$G(\mu) = \frac{x_1 x_3 + x_2 x_3 - 2x_1 x_2 + i\sqrt{3}x_3(x_2 - x_1)}{2x_3 - x_1 - x_2 + i\sqrt{3}(x_2 - x_1)}.$$

For $n = 4$, the situation is even simpler: $G(\mu)$ is the intersection between the hyperbolic geodesics $(x_1 x_3)$ and $(x_2 x_4)$.

Next we consider a Brownian flow $F_t$ in $M$, that is, a semimartingale flow with characteristic $(a(t, x, y, \omega), 0, t)$ as defined in Proposition 2.10, such that for every $x \in M$, $F_t(x)$ is a Brownian motion [which is equivalent to saying that $a(t, x, x, \omega) = \sum_{i=1}^d e_i \otimes e_i$, where $(e_i)_{1\leq i\leq d}$ is an orthonormal basis of $T_x M$]. We furthermore assume that $F_t$ is unstable, that is, whenever $y_1$ and $y_2$ are two distinct points in $M$, the distance between $F_t(y_1)$ and $F_t(y_2)$ converges to $+\infty$ as $t$ goes to infinity. Typical examples of unstable flows are given by isotropic Brownian flows in the hyperbolic plane with positive first Lyapounov exponent (see, e.g., [24]). The following result is immediate.

PROPOSITION 4.5. *Let $F_t$ be an unstable Brownian flow in $M$. Take $\varepsilon > 0$, $y_1 \neq y_2$ and let $X_t^i = F_t(y_i)$, $i = 1, 2$. Then for almost all $\omega$ there is $t(\omega)$ such that $\rho(X_t^1, X_t^2) \geq ((d-1)a - \varepsilon)t$ whenever $t \geq t(\omega)$.*

PROOF. For $x_1 \neq x_2$ and $u = (u_1, u_2) \in T_{(x_1,x_2)}(M \times M)$ write $u_i = v_i + w_i$, the orthogonal decomposition with $v_i$ tangential to the geodesic $(x_1 x_2)$. We have by the Hessian comparison theorem, $\nabla d\rho(u, u) \geq a \tanh(\frac{a\rho}{2})(\|w_0\|^2 + \|w_1\|^2)$ (see, e.g., [23] Lemma 1.1.1, where a similar calculation is done in positive curvature). As a consequence, the drift of $\Upsilon_t := \rho(X_t^1, X_t^2)$ is larger than $(d-1)a \tanh(\frac{a\Upsilon_t}{2})\, dt$. When $t$ is sufficiently large, then $\Upsilon_t$ is large by instability of the flow, so its drift is close to $(d-1)a\, dt$. On the other hand, since $\rho$ is 1-Lipschitz, the process $\Upsilon_t$ has a bracket satisfying $d\langle \Upsilon, \Upsilon \rangle_t \leq 2\, dt$. So we can conclude that for $\varepsilon > 0$ and $t$ large, $\Upsilon_t \geq ((d-1)a - \varepsilon)t$. □

PROPOSITION 4.6. *Let $X_t$ and $Y_t$ be two $M$-valued continuous functions converging to $X_\infty$ and $Y_\infty$ on $\partial M$, respectively, as $t \to \infty$. Suppose furthermore that $\lim_{t \to \infty} (\rho(o, X_t) - \ell t)/t = 0$, $\lim_{t \to \infty} (\rho(o, Y_t) - \ell t)/t = 0$ and $\rho(X_t, Y_t) \geq a't$ for $t$ large, for some constants $a', \ell > 0$ and $o \in M$. Then $X_\infty \neq Y_\infty$.*

PROOF. We only need to establish that if $\hat{P}_t$ is the orthogonal projection of $Y_t$ on the half-geodesic $[oX_t) := \{\varphi(o, X_t)(s),\, s \geq 0\}$, then $\lim_{t \to \infty} \rho(\hat{P}_t, Y_t) = \infty$. We will prove that, for $t$ large, $\rho(\hat{P}_t, Y_t) \geq \frac{a' \wedge \ell}{2} t$. When $\hat{P}_t = o$, this is clear. Let us consider the case where $\hat{P}_t \neq o$. Let $t$ be a time such that $\rho(o, \hat{P}_t) \geq \rho(o, X_t)$. Then

$$\rho(o, Y_t) \geq \rho(o, \hat{P}_t) = \rho(o, X_t) + \rho(X_t, \hat{P}_t) \geq \rho(o, X_t) + \rho(X_t, Y_t) - \rho(\hat{P}_t, Y_t),$$

which implies

$$\rho(\hat{P}_t, Y_t) \geq \rho(o, X_t) - \rho(o, Y_t) + \rho(X_t, Y_t).$$

This clearly implies $\rho(\hat{P}_t, Y_t) \geq \frac{a'}{2} t$ when $t$ is large. Now let $t$ be a time such that $\rho(o, \hat{P}_t) \leq \rho(o, X_t)$. We have

$$\rho(\hat{P}_t, Y_t) \geq \rho(o, Y_t) - \rho(o, \hat{P}_t),$$
$$\rho(\hat{P}_t, Y_t) \geq \rho(X_t, Y_t) - \rho(\hat{P}_t, X_t).$$

Adding the two together gives

$$2\rho(\hat{P}_t, Y_t) \geq \rho(o, Y_t) + \rho(X_t, Y_t) - \rho(o, X_t) \sim \rho(X_t, Y_t) \geq a't$$

which again proves that $\rho(\hat{P}_t, Y_t) \geq \frac{a'}{2} t$ when $t$ is large. Thus $X_\infty \neq Y_\infty$. □

A stochastic process $X_t$ is said to satisfy the law of large numbers with limit $\ell$ if $\lim_{t \to \infty} \frac{1}{t} \rho(o, X_t) = \ell$. It is known that all Brownian motions satisfy the law of large numbers with nonzero limit if (a) $M = \mathbb{H}^d$, (b) $M$ is the universal cover of a compact Riemannian manifold with negative curvature (see, e.g., [19]), or more generally, (c) $M$ is a Cartan–Hadamard manifold such that for some point $o \in M$, $\Delta \rho(o, -)$ converges to a negative constant as $r$ goes to infinity.

Note that if $\rho(X^i(t), y)$ are more or less all of the same length independent of $i$, that is, $\rho(X^i(t), o) = f(t) + R^i(t)$ with $R^i(t)$ small compared to $f(t)$, then the minimizer of $\sum_{i=1}^m \rho(X^i(t), y)$ is close to (or the same as) that of $\sum_{i=1}^m \rho^2(X^i(t), y)$. This is the consideration behind the following theorem.

THEOREM 4.7. *Let $\mu$ be a discrete probability measure on $M$ with finite support and no weight greater than or equal to $1/2$. Suppose $F_t$ is an unstable Brownian flow satisfying the conditions of Proposition 4.5 and such that $\lim_{t \to \infty} \frac{1}{t} \rho(o, X_t) = \ell$ where $\ell > 0$. Denote by $Z_t$ the exponential barycenter of the pushed forward measure $\mu_t \equiv F_t(\mu)$. Let $\mu_\infty$ be the measure $F_\infty(\mu)$ on $\partial M$. Then $\lim_{t \to \infty} Z_t = G(\mu_\infty)$ almost surely.*

PROOF. The measure $\mu_\infty$ is carried by a finite set, and by Propositions 4.5 and 4.6, for almost all $\omega$, $F_\infty(y_1, \omega) \neq F_\infty(y_2, \omega)$ if $y_1 \neq y_2$. This implies $\mu_\infty \in \mathcal{U}$. For $p \in M$ fixed denote by $\Theta_p(q)$ the point on $\partial M$ determined by the geodesic $\varphi(p, q)$, that is, $\Theta_p(q) = \phi(\dot{\varphi}(p, q)(0))$. Then $F_t(x)$ induces a measure $\bar{\mu}_{p,t}$ on $\partial M$ by $\Theta_p(\cdot)$ and $\bar{\mu}_{p,t} \to \mu_\infty$. Furthermore, by continuity of $G$, Proposition 4.3, $\lim_{t \to \infty} G(\bar{\mu}_{p,t}) = G(\mu_\infty)$.

Define
$$R_s = \int_M \left| \frac{1}{s} \rho(o, F_s(y)) - \ell \right| d\mu(y).$$

Since $\mu$ has finite support, $\lim_{s \to \infty} R_s = 0$ almost surely.

Set $Z'_t = G(\bar{\mu}_t)$. By definition $\int_{\partial M} \dot{\varphi}(Z'_t, x)(0) \, d\bar{\mu}_{o,t}(x) = 0$, so

$$\begin{aligned}
H_{\mu_s}(Z'_s) &= \int_M \dot{\gamma}(Z'_s, F_s(y))(0) \, d\mu(y) - \ell s \int_{\partial M} \dot{\varphi}(Z'_s, x)(0) \, d\bar{\mu}_{o,s}(x) \\
&= \int_M (\dot{\gamma}(Z'_s, F_s(y))(0) - \ell s \dot{\varphi}(Z'_s, F_s(y))(0)) \, d\mu(y) \\
&\quad + \ell s \int_M (\dot{\varphi}(Z'_s, F_s(y))(0) - \dot{\varphi}(Z'_s, \Theta_o(F_s(y)))(0)) \, d\mu(y) \\
&= \int_M (\rho(Z'_s, F_s(y)) - \ell s) \, \dot{\varphi}(Z'_s, F_s(y))(0) \, d\mu(y) \\
&\quad + \ell s \int_M (\dot{\varphi}(Z'_s, F_s(y))(0) - \dot{\varphi}(Z'_s, \Theta_o(F_s(y)))(0)) \, d\mu(y).
\end{aligned}$$

Set $R'_s = \ell \int_M \|\dot\varphi(Z'_s, F_s(y))(0) - \dot\varphi(Z'_s, \Theta_o(F_s(y)))(0)\| \, d\mu(y)$. By the convergence of $Z'_s$, the angular convergence of Brownian motions and the continuity off diagonal of the map $(z, x) \mapsto \dot\varphi(z, x)(0)$ from $M \times \bar{M}$ to $SM$, $\lim_{s \to \infty} R'_s = 0$ almost surely. We have

$$\|H_{\mu_s}(Z'_s)\| \leq \rho(o, Z'_s) + R_s s + s R'_s = \rho(o, Z'_s) + s(R_s + R'_s).$$

On the other hand, for $y \in M$, let $R_s(y) = |\frac{1}{s}\rho(o, F_s(y)) - \ell|$. Note that for $z \in B(Z'_s, 1)$,

$$\rho(z, F_s(y)) \geq \rho(o, F_s(y)) - \rho(z, o) \geq \ell s - R_s(y)s - [\rho(z, Z'_s) + \rho(Z'_s, o)]$$
$$\geq \ell s - R_s(y)s - 1 - \rho(Z'_s, o).$$

Consequently, by Lemma 4.1, if $u \in T_z M$ and $z \in B(Z'_s, 1)$,

$$\langle \nabla_u H_{\mu_s}, u \rangle$$
$$\leq -a\|u\|^2 \min_{w \in S_z M} \int_M \rho(z, F_s(y)) \sin^2(w, \dot\varphi(z, F_s(y))(0)) \, d\mu(y)$$
$$\leq -a\|u\|^2 \left( \ell s \min_{w \in S_z M} \int_M \sin^2(w, \dot\varphi(z, F_s(y))(0)) \, d\mu(y) \right.$$
$$\left. - R_s s - \rho(o, Z'_s) - 1 \right)$$
$$= -a\|u\|^2 (\ell s \alpha(z, \bar\mu_{z,s}) - R_s s - \rho(o, Z'_s) - 1),$$

where $\alpha$ is the continuous and positive function defined in Lemma 4.2. Since $\bar\mu_{z,s}$ converges to $F_\infty(\mu)$ and $Z'_s$ converges in $M$ to $G(\mu_\infty)$, there exist $C_1(\omega) > 0$ and $s(\omega) > 0$ such that for $s \geq s(\omega)$,

(4.8) $\qquad \langle \nabla_u H_{\mu_s}, u \rangle \leq -C_1 \|u\|^2 s \qquad \forall z \in B(Z'_s, 1), \ u \in T_z M.$

Next let $\beta_s(t)$, $t \geq 0$, be the $C^1$ path starting from $Z'_s$ and such that $\frac{\partial}{\partial t}\beta_s(t) = -H_{\mu_s}(\beta_s(t))$. Then with an argument similar to the end of the proof of Proposition 4.3, (4.8) implies that for every $t \geq 0$ smaller than the exit time from $B(Z'_s, 1)$,

$$\|H_{\mu_s}(\beta_s(t))\| \leq e^{-C_1 s} \|H_{\mu_s}(Z'_s)\|$$

and

(4.9)
$$\rho(\beta_s(t), \beta_s(0)) \leq \frac{\|H_{\mu_s}(Z'_s)\|}{C_1 s}$$
$$\leq \frac{\rho(o, Z'_s) + s(R_s + R'_s)}{C_1 s} = \frac{\rho(o, Z'_s)/s + R_s + R'_s}{C_1}.$$

Since the right-hand side converges to 0 as $s$ goes to infinity, we see that for $s$ large, $\beta_s(t)$ stays in $B(Z'_s, 1)$ for all $t \geq 0$. Moreover, $\beta_s(t)$ converges as $t$ goes to infinity to $\beta_s(\infty) = Z_s$ which is the only point in the manifold where $H_{\mu_s}$ vanishes. From (4.9), we get

$$\rho(Z_s, Z'_s) \leq \frac{\rho(o, Z'_s)/s + R_s + R'_s}{C_1}.$$

The right-hand side goes to 0, so $Z_s$ converges to $G(\mu_\infty)$. $\square$

REMARK 4.8. When $M$ is the hyperbolic space $\mathbb{H}^m$ and $a$ has second-order space derivatives almost surely bounded, Theorem 4.7 generalizes to a discrete measure $\mu$ with compact but not necessarily finite support. The proof is the same; the only difference is that one has to find another argument for the almost sure convergence of $R_s$ to 0. The new argument is as follows. In an exponential chart $\Psi$ based at $o$, let $(a', b', t)$ be the characteristic of the flow. Then $a'$ has second-order space derivatives almost surely bounded and $b' = \frac{1}{2}\Delta\Psi$ has first-order space derivatives almost surely bounded. Since the chart is centered at $o$, the distance to the origin is the same in the chart and in the manifold. Using [5], Theorem 2.1, we can say that almost surely, for all $y \in \text{supp}(\mu)$, $R_s(y) = |\frac{1}{s}\rho(o, F_s(y)) - \ell|$ is bounded by a constant depending only on the support of $\mu$. Since $R_s(y)$ converges almost surely to 0 as $s$ goes to infinity, by dominated convergence, $R_s$ goes to 0.

**Acknowledgments.** The authors would like to thank V. Kaimanovich, Y. LeJan, T. Lyons and M. Hairer for stimulating discussions and references. They would also like to thank the referees.

Département de mathématiques
Université de Poitiers
Téléport 2
BP 30179
F-86962 Futuroscope Chasseneuil Cedex
France
e-mail: marc.arnaudon@math.univ-poitiers.fr

Department of Computing
 and Mathematics
Nottingham Trent University
Nottingham NG1 4BU
United Kingdom
e-mail: xuemei.li@ntu.ac.uk